\documentclass[10pt]{amsart}
\usepackage{amsmath}
\usepackage{amsfonts}
\usepackage{amssymb}
\usepackage{graphicx}
\usepackage{amsthm,graphicx,color,yfonts}
\usepackage{pdfsync}
\usepackage{epstopdf}%

\usepackage[colorlinks=true]{hyperref}
\hypersetup{linkcolor=black,citecolor=black,filecolor=black,urlcolor=black} 

\providecommand{\U}[1]{\protect\rule{.1in}{.1in}}

\newcommand{\R} {\mathbb R}
\newcommand{\cuad}{{\sqcap\kern-.68em\sqcup}}

\newcommand{\ve}{\varepsilon}

\newcommand{\be}{\begin{equation}}
\newcommand{\ee}{\end{equation}}

\newcommand{\cal}{\mathcal}
\newcommand{\st}{{\tt t}}
\newcommand{\sy}{{\tt y}}
\newcommand{\pD}{{\mathring D_\tau}}

\definecolor{darkgreen}{rgb}{0.2,0.7,0.1}

\newcommand*{\ep}{\varepsilon}

\numberwithin{equation}{section}

\newtheorem{theorem}{Theorem}[section]
\newtheorem{proposition}{Proposition}[section]

\newtheorem{lemma}{Lemma}[section]

\newtheorem{remark}{Remark}[section]

\title[Rotationally symmetric  solutions to the Cahn-Hillard equation]{Rotationally symmetric  solutions to the Cahn-Hilliard equation}

\author{\'Alvaro Hern\'andez}
\address{Departamento de Ingenier\'{\i}a Matem\'atica and Centro de
Modelamiento Matem\'atico (UMI 2807 CNRS), Universidad de Chile, Casilla 170
Correo 3, Santiago, Chile.}
\email{ahernandez@dim.uchile.cl}
\thanks{\'A. H\'ernandez  was partially supported by Chilean research grants Fondecyt 109103 and Fondo Basal CMM-Chile, Project Anillo ACT-125}

\author{Micha{\l } Kowalczyk}
\address{Departamento de Ingenier\'{\i}a Matem\'atica and Centro
de Modelamiento Matem\'atico (UMI 2807 CNRS), Universidad de Chile, Casilla
170 Correo 3, Santiago, Chile.}
\email {kowalczy@dim.uchile.cl}
\thanks{M. Kowalczyk was partially supported by Chilean research grants Fondecyt 1090103, 1130126, Fondo Basal CMM-Chile, Project A\~nillo ACT-125 CAPDE}
\thanks{The Cahn-Hilliard equation and related to it the Allen-Cahn equation and the phase field model have been a subject of  extensive research of many mathematicians  for more than 30 years.  We have been a part of this group  and we owe it to Pauf Fife whose papers in the early 90ties  were for us an  introduction to the area and an inspiration for the present  work. For this reason we think it is appropriate to dedicate it to his memory}


\begin{document}

\maketitle

\begin{abstract}
This paper is devoted to construction of new solutions  to  the
Cahn-Hilliard equation in $\mathbb R^d$. Staring from a Delaunay unduloid
$D_\tau$ with parameter $\tau\in (0,\tau^*)$ we find for each sufficiently small
$\ve$ a solution $u$  of this equation which is periodic in the direction of the
$x_d$ axis and rotationally symmetric with respect to rotations about this axis.
The zero level set of $u$ approaches as $\ve\to 0$ the surface $D_\tau$. We use a refined version of the Lyapunov-Schmidt reduction method which simplifies very technical aspects of previous  constructions for similar problems. 
\end{abstract}

\section{Introduction} 
The Cahn-Hilliard  equation
\begin{equation}\label{CH evol}
\begin{aligned}
&u_t=-\Delta(\ep^2\Delta u-F'(u))\quad \text{ in }\Omega,\\
&\frac{\partial u}{\partial \nu}=0\quad \text{ on }\partial\Omega,\\
&\frac{\partial}{\partial \nu}(\ep^2\Delta u-F'(u))=0\quad \text{ on }\partial\Omega,
 \end{aligned}
\end{equation}
where $F$ is a {\it double-well potential}, is a model that describes the process of phase separation of two components of a binary  alloy. Here $\Omega\subset\R^{d}$, $d\geq 1$, is a bounded domain represents the place where the isolation of the components takes place, and $\nu$, as usual, denotes the outer normal on $\partial\Omega$. The function $u$ represents the concentration of one of the components and $\ep$ is the range of intermolecular forces. The double-well potential $F(u)$ corresponds to the free energy density at low temperatures, and in this paper we will take
\begin{equation*}
 F(u)=\frac14\left(1-u^2\right)^2,\quad F'(u)=u^3-u.
\end{equation*}
From now on we will denote $F'(u)=-f(u)$.

Equation \eqref{CH evol} can be derived from the gradient flow of the
Helmholtz-free energy functional 
\begin{equation}\label{helmholtz energy}
 E_{{\ep}}(u)=\int_\Omega \left(F(u(x))+\frac12\ep^2|\nabla
u(x)|^2\right)\ dx
\end{equation}
in $H^{-1}(\Omega)$ subject to the average concentration to be constant, i.e. 
\begin{equation}
\label{mass constr}
 \frac1{|\Omega|}\int_\Omega u\ dx =m, 
\end{equation} 
where $m\in[-1,1]$ (see \cite{MR1772733}, \cite{MR1901064} for details). Note that  constant functions $u\equiv \pm 1$ are minimizers of this functional subject to  $m=\pm 1$.

Stationary solutions of \eqref{CH evol}  satisfy the  Euler-Lagrange equation (with $f(u)=-F'(u)$)
\begin{equation}\label{CH stat}
\begin{aligned}
\ep^2\Delta u+f(u)=\delta_\ep\qquad &\text{ in }\Omega,\\
\frac{\partial u}{\partial \nu}=0\qquad &\text{ on }\partial\Omega,\\
\frac1{|\Omega|}\int_\Omega u\ dx=m
 \end{aligned}
\end{equation}
where $\delta_\ep$ is a Lagrange multiplier.

 Using $\Gamma$-convergence approach Modica \cite{MR866718}  showed that  minimizers $u_\ep$ of \eqref{helmholtz energy} subject to constraint (\ref{mass constr}) 
 $\Gamma$-converge, {as $\ep\to 0$}, to the function $1-2\chi_{A_0}$, where
$\chi_{A_0}$ is the characteristic function of an open set 
$A_0\subset\Omega$. Moreover $\partial A_0\cap \Omega$ is locally a surface of  constant mean curvature (CMC surface for short). Geometrically the set $A_0$ minimizes the perimeter functional $\mathrm{Per}_\Omega(A)$ among the sets $A\subset \Omega$ whose volume is fixed.  A generalisation of these results was given by Sternberg \cite{MR930124}.  Furthermore Hutchinson and Tonegawa \cite{MR1803974} studied limits of general critical points \eqref{helmholtz energy} and showed that their limits are locally minimal or CMC surfaces. On the other hand it  is known \cite{MR985990} that if a set $A\subset \Omega$  is an {\it isolated} mimimizer of the perimeter functional subject to the constant volume constraint then there exists a sequence of minimizers $u_\ep$ of (\ref{helmholtz energy}) which $\Gamma$ converges to $A$. This result can be used to construct solutions to (\ref{CH stat}) at least in dimension $2$ , see \cite{MR1399196}. The most complete construction is due to  Pacard and Ritor\'e \cite{MR2032110} 
who proved the following:  if $M$ is a  compact Riemannian manifold and $N$ is a non degenerate minimal or CMC sub manifold of $M$ which divides $M$ into $2$ disjoint components then for all sufficiently small  $\ve$ there exist critical points of \eqref{helmholtz energy} whose $0$ level set converges to $N$.  
The counterpart of this theory for the time dependent problem (\ref{CH evol}) was developed among others by Alikakos, Bates and Chen \cite{MR1308851} who proved that as $\ve\to 0$ the time evolution of interfaces is governed by the Helle-Shaw problem-of course CMC surfaces are stationary points of the flow. More detailed description of the Cahn-Hilliard flow and key spectral tools   can be found for instance in  \cite{MR1237062}, \cite{MR1382059}, \cite{MR1613496}, \cite{MR1797871}, \cite{MR1284813} and the references therein. Additional examples of stationary solutions for the singular perturbation problem  in a bounded domain have been constructed  in \cite{MR1632937}, \cite{MR1636688}, \cite{MR1737000}.

Scaling variables $x\mapsto x/\ep$ in  (\ref{CH stat})  and letting $\ep$ tend to $0$ leads in a natural way to the following problem:
\begin{equation}
\label{CH1}
\Delta u+f(u)=\delta, \quad \mbox{in}\ \R^d.
\end{equation}
In dimension $d=1$ there is an  obvious solution of this problem when $\delta=0$, namely the unique odd and monotonically increasing heteroclinic solution $H$ of the ODE, which satisfies
\begin{equation}
\begin{aligned}
H''+f(H)&=0, \quad \mbox{in}\ \R, \\
H(\pm \infty)&=\pm 1.
\end{aligned}
\label{hetero}
\end{equation} 
If ${\tt a}\in \R^d$ is a unit vector and $b\in \R$ then the function
\[
u(x)=H({\tt a}\cdot x+b),
\] 
is also a solution of (\ref{CH1}) with $\delta=0$. When $\delta\neq 0$ there exist radially symmetric solutions to (\ref{CH1}), see \cite{MR682268}. In both cases the level sets of the solutions are CMC surfaces, in the former case their mean curvature  is $0$ and in the latter case it is positive number equal to $\frac{1}{d-1}\frac{1}{R_0}$, where $R_0$ is the radius of the level set of the solution.  The radially symmetric solutions in $\R^{d-1}$ can be lifted trivially to $\R^d$ giving solutions whose nodal sets are cylinders, which again are CMC surfaces.

Dilating of the independent  variable by a (large) factor $\ep^{-1}>0$
\[
{\tt x}\longmapsto \ep^{-1}{\tt x},
\]
we obtain an equivalent form of (\ref {CH1}):
\begin{equation}
\label{CH}
\ep\Delta u+\frac{1}{\ep}f(u)={{\ell}}_\ep, \quad \mbox{in}\
\R^{{^d}}.
\end{equation} 
where we have denoted $\frac{\delta}{\ep}={{\ell}}_\ep$. Clearly, if
$u_\ep$ is a solution of (\ref{CH}) then $v({\tt x})= u_\ep\big({\ep}{\tt
x}\big)$ is a solution of (\ref{CH1}). On the other hand, if $v$ is a solution
of (\ref{CH1}) then $u_\ep({\tt x})=v\big(\frac{\tt x}{\ep}\big)$ is a solution
of (\ref{CH}). In particular this means that while phase transition of the
solutions of (\ref{CH1}) are of order $1$, for the solutions of (\ref{CH}) they
are of order $\ep$. Thus the latter are more ``concentrated". In the sequel we
will focus on solving (\ref{CH}).  From what we have said above about the
singular perturbation problem  it is clear that level sets of these solutions
should converge, as $\ep$ tends to $0$, to CMC surfaces in $\R^3$. In fact we
expect (on the basis of formal calculations in section \ref{sec formal}) that
the Lagrange multiplier
\[
{{\ell}}_\ep=-\frac{1}{2}H_\varSigma\int_\R H'(s)^2\,ds+\mathcal O(\ep),
\]
where $\varSigma$ is the surface of the phase transition and $H_\varSigma$ is its mean curvature.

We will now introduce a family of embedded CMC surfaces which are good candidates to be the limits of the nodal surfaces.
We recall that  Delaunay {\it unduloids}  \cite{delaunay_1},  \cite{MR869541} are a one parameter family $D_\tau$, $\tau\in (0,1)$ of embedded, periodic CMC surfaces of revolution in $\R^3$. When the real parameter $\tau$ tends to $1$ the surfaces $D_\tau$ approach the straight cylinder while  when $\tau\to 0$ they become an array of identical spheres arranged along the $x_3$ axis.
It turns out that  Delaunay  surfaces  can be constructed in  any dimension $d>3$ and from now on by $D_\tau$, $\tau\in (0, \tau_*)$ we will denote the family of Delaunay surfaces in $\R^d$. We note that the parameter $\tau_*$ satisfies:
\begin{equation}
\label{tau star}
\tau_*=\frac{2}{(d-1)(d-2)^{(d-2)}}.
\end{equation}
Again, in the limit $\tau\to \tau_*$ the surfaces $D_\tau$ approach the straight cylinder, see \cite{MR2590386} or Section \ref{sec delaunay} for details. 
%

It is convenient to "normalize'' the Delaunay surface and suppose that the mean
curvature of $D_\tau$ is $1$ for all $\tau\in (0, \tau_*)$.  We will also denote
by  $N_\tau$  the  vector field normal  to  $D_{\tau}$.
Let us notice that the surface $D_{\tau}$ divides the space into two disjoint  components $\Omega^\pm_\tau$, such that  $\R^d\setminus D_{\tau}={\Omega}_{\tau}^+\cup{\Omega}_{\tau}^-$, where  $N_\tau$ points towards ${\Omega}_{\tau}^+$. By changing the orientation of $D_\tau$ if necessary  we can  chose $N_\tau$ in such a way that $\Omega_\tau^+$ contains the $x_d$ axis. 

Our  result is:
\begin{theorem}\label{teorema 1}
 For all $\tau\in(0,\tau_*)$ there exits $\ep_\tau>0$ such that for all $\ep\in (0, \ep_\tau)$  the problem 
\begin{equation}\label{CH ep}
\ep \Delta u+\frac{1}{\ep}f(u)={{\ell}}_\ep\quad\text{ in}\ \R^d
\end{equation} 
has a solution $u_{\tau, \ep}$, which is  one-periodic in the direction of the 
$x_d$-axis and rotationally symmetric with respect to rotations about the same
axis. As $\ep\to 0$ we have  ${{\ell}}_\ep =1+\mathcal O(\ep)$, and 
$u_{\tau, \ve}$  satisfies
\begin{align*}
u_{\tau, \ep}\to 1 \text{ as }\ep\to 0&\text{ in }{\Omega}_{\tau}^+,\\
u_{\tau, \ep}\to -1\text{ as }\ep\to 0&\text{ in }{\Omega}_{\tau}^-,
\end{align*} 
uniformly over compacts. 
\end{theorem}
\begin{remark}
In this paper we took $f(u)=u-u^3$, which is the standard  nonlinearity for the Cahn-Hilliard equation. Theorem \ref{teorema 1} holds  for more general nonlinearities of bistable, balanced type, namely $f\in C^3$ such that $f(u)=-F'(u)$ where $F$ is a double well, even potential with non degenerate wells at $\pm 1$.  Rather straightforward modifications required in the proof of the more general setting are left to the reader.
\end{remark}

We will explain  now the implementation  of the Lyapunov-Schmidt reduction we used in this paper and discuss the differences between our approach and the older implementations  which can be found in \cite{MR2032110} and \cite{dkp_dg}, \cite{costa}. Let us first recall the standard Lyapunov-Schmidt reduction method in its abstract version (see \cite{MR660633}). Given  Banach spaces $X, Y$ and a linear operator $A\colon X\to Y$ and a continuous, nonlinear operator $N\colon X\to Z$,  we are to solve the problem:
\begin{equation}
\label{lyap 1}
Ax-N(x)=0. 
\end{equation}
Let  
\[
\mathcal N(A)=Y\subset X, \qquad \mathcal R(A)=W\subset Z,
\]
and let $\pi_Y$, $\pi_W$ be the projections on the corresponding subspaces.  There exists a bounded linear operator $K\colon W\to \mathcal R(I-\pi_Y)$ (the right inverse of 
$A$) such that $AK=I$ on $W$ and $KA=I-\pi_Y$, and moreover the equation (\ref{lyap 1}) is equivalent to the equation
\begin{equation}\label{lyap 2}
\begin{aligned}
&x=y+z, \qquad y\in Y,\quad z\in \mathcal R(I-\pi_Y)\\ 
&\qquad z-K\pi_W N(y+z)=0, \\
&\qquad (I-\pi_W) N(y+z)=0.
\end{aligned}
\end{equation}
In applications the Lyapunov-Schmidt method consists of reducing (\ref{lyap 1}) to (\ref{lyap 2}), solving the first equation for $z$ with $y$ given (which usually can be done by a fixed point argument) and  replacing this solution in the second equation to obtain {\it the reduced problem}
\begin{equation}
\label{lyap 3}
(I-\pi_W) N(y+z(y))=0.
\end{equation}
In practice several complications may arise and we will illustrate this considering a related to our problem which was treated by Pacard and Ritor\'e  \cite{MR2032110}, and in many aspects it is similar to problem we consider in this paper. Let $M$ be a compact, closed manifold  of dimension $n$ and $N\subset  M$ a minimal $n-1$ dimensional  sub manifold  which divides $M$ into two disjoint components. Consider the problem
\begin{equation}
\label{lyap 4}
\ve^2\Delta_M u+u(1-u^2)=0, \qquad \mbox{on}\ M.
\end{equation} 
We say that   $N$ is non degenerate if the Jacobi operator of $N$ 
\[
J_N=\Delta_N+|A_N|^2+\mathrm{Ric}_g\,(\nu_N, \nu_N)
\]
has empty kernel ($\Delta_N$ is the Laplace-Beltrami operator on $N$, $|A_N|^2$ is the norm of the second fundamental form, $\mathrm{Ric}_g$ is the Ricci tensor on $M$ and $\nu_N$ is the normal vector to $N$). The result proven in \cite{MR2032110} is: given a non degenerate, minimal sub manifold $N$ of $M$ for each sufficiently small $\ve$ there exists a solution $u_\ve$ of (\ref{lyap 4}) such that the zero level set of $u_\ve$ approaches $N$ as $\ve\to 0$. Moreover, $u_\ve$ converges to $\pm 1$ uniformly over compacts of the two disjoint components of $M\setminus N$.

Let us explain now the  implementation of the Lyapunov-Schmidt reduction  in  \cite{MR2032110}. It is expected that for $x\in M$ near $N$ we should have $u_\ve(x)=H(\ve^{-1}\mathrm{dist}\,(x,N))+\varphi$, where $\mathrm{dist}\,(\cdot,N)$ is the signed distance function on $M$, $H$ is the unique odd, monotonically increasing solution of $-H''=H(1-H^2)$ in $\mathbb R$ and $\varphi$ is a small perturbation. The problem to solve for $\varphi$ amounts to inverting the linearized operator  around $H(\ve^{-1}\mathrm{dist}\,(x,N))$ which has form
\[
L=\Delta_M +f'(H(\ve^{-1}\mathrm{dist}\,(\cdot,N))).
\]
It is known that the norm of $L^{-1}$ is large due to local translational invariance of the problem. Thus we need to perturb $N$ as well. To describe this perturbation we consider a manifold $N_h$ to be a normal graph over $N$ described by a smooth and small function $h\colon N\to \R$. Furthermoe we let $t_h(x)=\mathrm{dist}\, (x, N_h)$ to be the signed distance from $N_h$. Then  we look for a solution of the form 
\[
u= H\Big(\frac{t_h}{\ve}\Big) +\varphi.
\]
Now both $h$ and $\varphi$ are unknowns. The problem to solve for $\varphi$ is 
\[
L_h\varphi=\mathcal F(h,\varphi),
\]
where $L_h$ is the linearized operator around $H\Big(\frac{t_h}{\ve}\Big)$. The Lyapunov-Schmidt reduction strategy amounts to projection of the above equation onto the function $H'\Big(\frac{t_h}{\ve}\Big)$ and its complement, denote this last projection by $\pi_h$. This leads to a problem for $\varphi$
\[
\pi_h L_h\varphi=\pi_h \mathcal F(h,\varphi),
\]
which we solve first for a given $h$,  
and the problem for $h$
\begin{equation}
\label{lyap 5}
J_{N_h} h=\mathcal G(h),
\end{equation}
which we solve next ($J_{N_h}$ is the Jacobi operator of $N_h$). Let us discuss (\ref{lyap 5}). We notice that the expression of $J_{N_h}$ in local coordinates will depend in general on  $h$ and its derivatives up to order $3$, while the Jacobi operator is itself only a second order operator. This loss of regularity was dealt with in  \cite{MR2032110} using a regularisation procedure. In a series of papers \cite{dkp_dg}, \cite{costa}, \cite{nested} del Pino, Kowalczyk and Wei introduced a slightly different approach to circumvent this problem.  It amounts to considering perturbation in the normal direction of the {\it fixed} manifold $N$ so that $u=H\Big(\frac{t +h}{\ve}\Big)+\dots$, where now $t$ is the signed distance from $N$ and $h$ is a smooth, unknown function defined on $N$. Equation (\ref{lyap 5}) takes form 
\begin{equation}
\label{lyap 6}
J_{N} h=\mathcal G(h),
\end{equation}
and the problem of the loss of regularity is thus avoided. The problem is now reduced to finding a fixed point of $J_N\circ \mathcal G (h)$, using for example   Banach fixed point theorem. To do this we need to know that $\mathcal G$ is at least Lipschitz in $h$. In both implementations of the Lyapunov-Schidt reduction described above this is rather complicated technical point since $\mathcal G$ depends in a non explicit, non local and non linear way on $h$. This is mainly due to the fact that the linearized operator $L_h$ still depends on $h$ through the potential $f'(H(\frac{t +h}{\ve}))$. 

In this paper we propose still another modification to the method. The idea is simple: instead of working with an approximation of the form $u=H\Big(\frac{t +h}{\ve}\Big)+\dots$ with $h$ unknown we will improve the initial approximation to $w(t, y)=H\big(\frac{t}{\ve}\big)+\dots$, $t$ being the signed distance to $N$ and $y\in N$ in such a way that we do not need to "move" $N$ anymore. In other words $h$ will be determined with some sufficient precision before setting up the Lyapunov-Schmidt reduction, which with this modification will  look  like the abstract setting described at the beginning. This way we avoid both the loss of regularity and technical difficulties due to complicated character of the nonlinear function $\mathcal G(h)$.
This is described in detail in section \ref{ls}.

\section{Preliminaries}
\subsection{The surfaces of Delaunay}\label{sec delaunay}\setcounter{equation}{0} 

The Delaunay unduloids $D_\tau$, $\tau\in (0, \tau_*)$   are CMC surfaces of revolution in $\R^d$. Thus for instance in $\R^3$ one can  parametrize them in the form 
\[
{\tt x} (t, \theta)=(\rho(t)\cos\theta, \rho(t)\sin\theta, t), 
\]
where $\rho(t)$ solves
\[
\rho_{tt}-\frac{1}{\rho}(1+\rho_t^2)-(1+\rho_t^2)^{3/2}=0.
\]
However, in this paper we will use mostly isothermal coordinates of $D_\tau\subset \R^d$:
\begin{equation}
\label{dela 1}
X_\tau(s, \varTheta)=\frac{1}{2}(\tau e^{\,\sigma_\tau(s)}\varTheta, \kappa_\tau(s)), \quad (s, \varTheta)\in \R\times S^{d-2},
\end{equation}
where functions $(\sigma_\tau, \kappa_\tau)$ are the unique solutions of the following system of ODEs:
\begin{equation}
\label{dela 2}
\begin{aligned}
(\partial_s\sigma_\tau)^2+\frac{1}{4}\tau^2(e^{\,\sigma_\tau}+e^{\,
(2-d)\sigma_\tau})^2=0, &\quad \partial_s\sigma_\tau(0)={0}, \quad
\sigma_\tau(0)<0, \\
\partial_s\kappa_\tau-\frac{1}{4}\tau^2(e^{\,\sigma_\tau}+e^{\,(2-d)\sigma_\tau})=0, &\quad \kappa_\tau(0)=0.
\end{aligned}
\end{equation}
We will now summarize some basic facts about the Delaunay surfaces and their isothermal parametrization. The function $\sigma_\tau$ is periodic, and consequently the surfaces $D_\tau$ are one-periodic along the $x_d$-axis: namely if $2T_\tau$ denotes the minimal  period then
\[
D_\tau=D_\tau+2T_\tau {\tt e}_d.
\]
Clearly we have the relation 
\[
T_\tau=\frac{1}{4} \kappa_\tau(2s_\tau),
\]
where $2s_\tau$ is the minimal period of $\sigma_\tau$. 

The Jacobi operator $\mathcal J_\tau$ of $D_\tau$ is defined by:
\begin{equation}
\mathcal J_\tau:=\Delta_{D_\tau}+|A_\tau|^2,
\label{def jacobi}
\end{equation}
where $\Delta_{D_\tau}$ is the Laplace-Beltrami operator on $D_\tau$ and  $|A_\tau|^2$ is the square of the norm of the second fundamental form of $D_\tau$. The Jacobi operator is of fundamental importance in this paper and to understand well its properties we will first consider the special case $d=3$. In the isothermal coordinates $(s, \theta)\in \R\times S^1$ its expression is given by:
\begin{equation}
\label{jacob 1}
L_\tau=\frac{1}{\tau^2 e^{\,2\sigma_\tau}}\left\{\partial_s^2+\partial_\theta^2+\tau^2\cosh (2\sigma_\tau)\right\}.
\end{equation}
The Jacobi fields on $D_\tau$, which are elements of the kernel of
${L}_\tau$ are of three types:

\begin{itemize}
\item[(1)] {\it Jacobi fields arising from  infinitesimal translations.} For any
 ${\tt e}\in \R^3$, $|{\tt e}|=1$ {
the constant Killing field associated to
translations   
\begin{equation*} 
{\tt x}\longmapsto {\tt e}
\end{equation*}
induces the following Jacobi fields}

\[
\Phi_\tau^{T, {\tt e}}={\tt e}\cdot N_\tau,
\]
where $N_\tau$ is the unit normal vector to $D_\tau$. The coordinate vectors ${\tt e}_j$, $j=1,2,3$ generate three linearly independent Jacobi fields $\Phi^{T, {\tt e}_j}_\tau$ corresponding to translations of $D_\tau$ in the directions of the coordinate axis. 
We note that in the isothermal coordinates  
\[
\Phi_\tau^{T, {\tt e}_3}=\Phi_\tau^{T, {\tt e}_3}(s), \quad \Phi_\tau^{T, {\tt e}_j}=\Phi_\tau^{T, {\tt e}_j}(s, \theta), \quad j=1,2.
\]
It is important to notice that the Jacobi fields  $\Phi_\tau^{T, {\tt e}_j}$ are bounded.
\item[(2)]{\it Jacobi fields arising from infinitesimal  rotations.} 
Let ${\tt e}\in \R^3$, $|{\tt e}|=1$ be such that ${\tt e}\cdot {\tt e}_3=0$. The Killing vector field corresponding to the rotation about the vector ${\tt e}$ is:
\[
{\tt x}\longmapsto ({\tt x}\cdot{\tt e}){\tt e}_3-({\tt x}\cdot{\tt e}_3){\tt e}.
\] 
We define the Jacobi field associated to this vector field by:
\[
\Phi_\tau^{R, {\tt e}}=\left[({\tt x}\cdot{\tt e}){\tt e}_3-({\tt x}\cdot{\tt e}_3){\tt e}\right]\cdot N_\tau.
\]
There are clearly two linearly independent Jacobi fields associated to the rotations. They are:
\[
\Phi_\tau^{R, {\tt e}_1}, \quad \Phi_\tau^{R, {\tt e}_2},
\]
and they correspond to  rotations about the coordinate axis. Note that in isothermal coordinates functions $\Phi_\tau^{R, {\tt e}_j}$, $j=1,2$ grow linearly as functions of $s$.  
\item[(3)]{\it Jacobi field associated with the variation of the Delaunay parameter.}
We define:
\[
\Phi^D_\tau=-\partial_\tau X_\tau\cdot N_\tau.
\]
This Jacobi field is somewhat harder to write explicitly however it can be  shown that the function $\Phi^D_\tau(s)$ is linearly growing.
\end{itemize}
In summary, the Jacobi operator ${L}_\tau$ has at least $6$ explicit
Jacobi fields which are either linearly growing or bounded. By a result of 
Mazzeo and Pacard \cite{MR1941630} we know that these are {\it all} Jacobi
fields with at most linear growth. To explain this let us observe that by
separation of variables the equation $J_\tau\varphi =0$ separates into a
sequence of problems
\[
{L}_{\tau, j}\varphi=0, \quad {L}_{\tau, j}=\frac{1}{\tau^2
e^{\,2\sigma_\tau}}\left\{\partial_s^2+\tau^2\cosh (2\sigma_\tau)-j^2\right\},
\quad |j|=0, 1, \dots.
\]
Then we have:
\begin{proposition}[\cite{MR1941630}]\label{proposition 2.1}
The homogeneous problem ${L}_{\tau, j}\varphi=0$ has the following
solutions:
\begin{enumerate}
\item
one periodic and one linearly growing solution when $j=0$ or $|j|=1$;
\item
two solutions $\varphi^\pm_{\tau, j}(s)$ which satisfy:
\[
\varphi^\pm_{\tau, j}(s+s_\tau)=e^{\,\pm \zeta_{\tau, j} s_\tau} \varphi^\pm_{\tau, j}(s),
\]
with 
\[
\gamma_{\tau, j}=\mathrm{Re}\,\zeta_{\tau, j}>0,
\]
when $|j|>2$.
The numbers $\zeta_{\tau, j}$ are  the indicial roots of the operators $L_{\tau, j}$ and correspond to the behavior of the solutions of the homogeneous problem at $\pm \infty$.  
\end{enumerate}
\end{proposition} 
This basic facts can be generalized for  the Jacobi operator of Delaunay surfaces in $\R^d$, $d>3$. We will summarize them now and  refer the reader to \cite{MR2590386} for details. We have the following at most linearly growing Jacobi fields:
\begin{itemize}
\item[(1)]  The are $d$ bounded, periodic Jacobi fields arising from 
infinitesimal translations. They will be  denoted by $\Phi_\tau^{T, {\tt e}_j}$,
$j=1,\dots, d$.
\item[(2)] There are $d-1$ Jacobi fields arising from infinitesimal  rotations
in the direction of  $x_j$ axis, $j=1, \dots, d-1$. We will denote them by 
$\Phi_\tau^{R, {\tt e}_j}$, $j=1, \dots, d-1$. 
and they correspond to  rotations about the coordinate axis. Note that in isothermal coordinates functions $\Phi_\tau^{R, {\tt e}_j}$,  grow linearly as functions of $s$.  
\item[(3)]There is one Jacobi field associated with the variation of the
Delaunay parameter $\Phi^D_\tau=-\partial_\tau X_\tau\cdot N_\tau$, and it is
linearly growing.
\end{itemize}

\subsection{Fermi  coordinates and shifted Fermi coordinates near  a CMC surface}\label{sec fermi}

\setcounter{equation}{0} 

Let $\varSigma$ be an embedded CMC surface in $\R^d$ and let  $H_\varSigma$ denote its mean curvature.  By $N$ we will denote its unit normal. We will assume that there exists  a tubular neighborhood $\mathcal N_\delta$   of $\varSigma$ of width  $2\delta$ in which we can introduce local system of coordinates (Fermi coordinates) $(y, z)\in \varSigma\times (-\delta, \delta)$ by setting:
\[
{\tt x}\longmapsto (y, z), \quad \mbox{where}\ {\tt x}=y+z N(y).
\]
We suppose that this map, which we denote by $Y$, is a diffeomorphism from $\mathcal N_\delta$ to $\varSigma\times (-\delta, \delta)$ whenever $\delta$ is taken sufficiently small. In the sequel we will use the inverse of this  map
\[
\begin{aligned}
Y^{-1}\colon\varSigma\times(-\delta, \delta)&\longrightarrow \mathcal N_\delta\\
(y, z)&\longmapsto {\tt x}.
\end{aligned}
\]
Given a function $w\colon\mathcal N_\delta\to \R^d$ we define its pullback $Y^*w$ to $ \varSigma\times(-\delta, \delta)$ by the diffeomorphism $Y$ as:
\[
Y^* w({y,z})=w\circ Y^{-1}(y,z).
\]
For technical reasons we will chose later the size of the tubular neighbourhood $\delta$ depending on $\ve$ but for now on we just take $\delta$ small. 

Next we will define  {\it shifted} Fermi coordinates. To do this we let  $h\colon \varSigma\to \R$ be a given smooth function such that the map
\[
{\tt x}\longmapsto (y, t), \quad\mbox{where}\ {{\tt x}}= 
y+\left(t+h(y)\right) N(y),
\]
is a diffeomorphism from $\mathcal N_\delta$ into $\varSigma\times(-\delta, \delta)$. We will denote this map by  $Y_h$ and by  $Y^{-1}_h$  we will denote its inverse, finally  by $Y^*_hw$ we will denote the pullback of $w\colon\mathcal N_\delta\to \R^d$ by $Y_h$:
\[
Y^*_h w(y,t)=w\circ Y_h^{-1}(y, t).
\]  
 
It will be convenient to have at hand expressions for the Laplacian  in Fermi and shifted Fermi  coordinates. To derive them by $\varSigma_z$ we will denote the surface $\varSigma+ zN$ i.e. the original surface $\varSigma$ shifted in the direction of the normal by $z$. Locally near $\varSigma$ we have
\[
\Delta = \Delta_{\varSigma_z}+\partial_z^2- H_{\varSigma_z}\partial_z.
\]
We denote $z{\mathbb B}_{\varSigma,
z}=\Delta_{\varSigma_z}-\Delta_\varSigma$. The operator  ${\mathbb
B}_{\varSigma, z}$ is a second order differential operator. To expand the
curvature term we use the well known formula:
\[
H_{\varSigma_z}=\sum_{j=1}^{d-1} \frac{{\tt k}_j}{1-z{\tt
k}_j}=H_\varSigma+z|A_\varSigma|^2+z^2 \sum_{j=1}^{d-1} {\tt
k}_j^3+{{\cal O}(z^3)}=H_\varSigma+z|A_\varSigma|^2+z^2{\mathbb
Q_{\varSigma, z}}.
\]
where ${\tt k}_j$ are the principal curvatures of $\varSigma$.  In summary we have:
\[
\Delta=\Delta_\varSigma{+\partial_z^2-(H_\varSigma +
z|A_\varSigma|^2)}\partial_z+z{\mathbb B}_{\varSigma, z}+z^2\mathbb
Q_{\varSigma, z}.
\]
The reason we expanded the Laplacian in this way will become clear later on. From this it is easy to obtain a formula for the Laplacian in the  shifted Fermi coordinates:
\begin{equation}
\label{fermi 1}
\Delta=\Delta_\varSigma+ (1+|\nabla_\varSigma h|^2)\partial_t^2-(H_\varSigma+\Delta_\varSigma
h+(t+h)|A_\varSigma|^2)\partial_t+(t+h)\mathbb B_{\varSigma, t+h}+(t+h)^2\mathbb
Q_{\varSigma, t+h}.
%
\end{equation}

Anticipating the content of the next section we introduce the stretched shifted Fermi coordinate 
\[
{\tt t}=\frac{t}{\ep}, \quad {\tt y}=y.
\]
Formal consideration will show that an approximate solution $w_\ve$ of the Cahn-Hilliard can be obtained if we assume that it is a function of the form:
\[
Y^*_h w_\ve(y, t)= H\left(\frac{t}{\ve}\right)+o(1),
\] 
where $H$ is the heteroclinic solution of \eqref{hetero}

As before we have a diffeomorphism $Y_{\ve, h}$ and its inverse  $Y^{-1}_{\ve, h}\colon\varSigma\times(-\frac{\delta}{\ve}, \frac{\delta}{\ve})\to \mathcal N_{\delta}$, and for any function $w\colon\mathcal N_{\delta}\to \R^k$ we define its pullback by $Y_{\ve, h}$ by:
\[
Y^*_{\ve, h} w({\tt y}, {\tt t})=w\circ Y^{-1}_{\ve, h}({\tt y}, {\tt t}).
\] 
Taking onto account formula (\ref{fermi 1}) we get
\begin{equation}
\label{fermi 2}
\Delta=\Delta_\varSigma+\ve^{-2}(1+|\nabla_\varSigma h|^2)\partial_{\tt
t}^2-\ve^{-1}(H_\varSigma+\Delta_\varSigma h+(\ve{\tt
t}+h)|A_\varSigma|^2)\partial_{\tt t}\\+(\ve{\tt t}+h) \mathbb B_{\varSigma, \ve
{\tt t}+h}+(\ve{\tt t}+h)^2\mathbb Q_{\varSigma, \ve {\tt t}+h}.
\end{equation}

\subsection{Formal expansion of the solution of the Cahn-Hilliard equation concentrating on $\varSigma$}\label{sec formal}
For the purpose of formal calculations we will assume that the the solution of ({\ref{CH}) near $\varSigma$ is a function $w$, which depends  on the stretched and shifted Fermi coordinates $({y}, {\tt t})$, in the following way
\begin{equation}
Y^*_{\ve, h} w({\tt y}, {\tt t})=U({\tt t})+\ve^2\psi_0(\sy,\st),,
\label{ansatz}
\end{equation}
for some functions $U$ and $\psi_0$  which we will determine. 
Moreover, we will assume that 
\[
h=\ve^2 h_0,
\]
where $h_0$ is a constant to be chosen. 

{ 
To determine $U$ and $\psi_0$ we write the error 
$N_\ve(w)-\ell_\ve:=\ve\Delta w+\frac1\ve w(1-w^2)-\ell_\ve$ in the local coordinates
}
{
\begin{equation}\label{Error}
\begin{split}
N_\ve(w)-\ell_\ve&=\ve\big\{\ve^{-2} \partial_{\tt
t}^2U-\ve^{-1}H_\varSigma\partial_{\tt
t}U+\ve^{-2}f(U)-\ve^{-1}\ell_\ve\\
&+\partial_{\tt t }^2\psi_0- \ve H_\varSigma\partial_{\tt
t}\psi_0+f'(U)\psi_0-({\tt t}+\ve h_0)|A_\varSigma|^2\partial_{\tt t}U\\
&+\frac1\ve f(w)-\frac1\ve f(U)-\ve f'(U)\psi_0\\
&+\ve^2\Delta_\varSigma\psi_0-\ve^2( {\tt t}+\ve
h_0)|A_\varSigma|^2\partial_{\tt t}\psi_0)\\
&+\left((\ve{\tt t}+h){\mathbb B}_{\varSigma, \ve{\tt t}+h}+(\ve{\tt t}+h)^2{\mathbb
Q}_{\varSigma,\ve{\tt t}+h}\right)(U+\ve^2\psi_0)\big\}.
\end{split}
\end{equation}
}
{
In order to get as small as possible this approximation we have to get rid of
the first three terms  of the right hand side of the expression above, since they show the lower powers
in $\ve$. To write things compactly let:  
\begin{align*}
S_0(w)&:=\partial_{\tt t}^2 w-\ve H_\varSigma\partial_{\tt t}w+f(w), \\
L_0 w&:=\partial_{\tt t}^2 w-\ve H_\varSigma\partial_{\tt t}w+f'(U)w.
\end{align*}
}
{
With this notation and the ansatz \eqref{ansatz} we can write the problem in the
form:
\[
S_0(U)+\ve^2 L_0\psi_0+ Q_0(U+\ve^2\psi_0)=\ve\ell_\ve.
\] 
where $Q(U+\ve^2\psi_0)$ represents the rest of the terms in \eqref{Error}, and
we also have to determine the  Lagrange multiplier $\ell_\ve$.  
}

Thus we take the Lagrange multiplier $\ell_\ve$ to
be a number such that the following ODE has a  unique, monotonically increasing
solution $U$
\begin{equation}
\label{eq U}
\begin{aligned}
{S_0(U)=}U''-\ve H_\varSigma U'+f(U)&=\ve\ell_\ve, \quad \mbox{in}\ \R,\\
f(U(\pm \infty))&=\ve\ell_\ve.
\end{aligned}
\end{equation}
This function is easy to find  by
perturbing the heteroclinic solution $H$. 

Note that since   $f$ is odd symmetric we have $U(\pm\infty)=\pm 1+\sigma_\ve$,
where 
\[
f(\pm 1+\sigma_\ve)=\ve\ell_\ve.
\]
Also, we have
\[
\ell_\ve=\ell_0+\mathcal O(\ve), \quad \ell_0=-\frac{1}{2}H_\varSigma\int_\R H'(s)^2\,ds,
\]
and 
\[
U({\tt t})=H({\tt t}) +\mathcal O(\ve).
\]

Next, we will determine the $\mathcal O(\ve^2)$ correction $\psi_0$.
Ignoring terms of order ${{\cal O}}(\ve^3)$ we get the following equation
to solve:
\begin{equation}
\partial_{\tt t}^2 \psi_0-\ve H_\varSigma\partial_\st \psi_0+(1-3U^2)\psi_0=(\st +\ve h_0)|A_\varSigma|^2 \partial_{\tt t} U.
\label{match 3}
\end{equation}
It is convenient to consider, more generally, an ODE (with the right hand side
possibly depending on the variable $\sy$) of the form:
\begin{equation}
\label{match 5}
\partial_{\tt t}^2 \varphi-\ve H_\varSigma\partial_\st \varphi+(1-3U^2)\varphi=g(\sy, {\tt t}).
\end{equation}
A solution of this problem  can be found by the variation of parameters formula. Indeed, the fundamental set of the ODE is spanned by the functions 
\begin{equation}
\label{match 6}
  \partial_\st U=\mathcal O\big((\cosh \st)^{\eta^\pm}\big),
\quad W(\st)=\mathcal O\big((\cosh \st)^{\nu^\pm}\big),\quad \st\to \pm \infty.
\end{equation}
with 
\begin{equation}\label{indicial 1}
\eta^\pm=\frac{1}{2}\Big(\ve H_\varSigma-\sqrt{-4\iota(\pm\infty)+\ve^2 H^2_\varSigma}\Big), \quad
\nu^\pm=\frac{1}{2}\Big(\ve H_\varSigma+\sqrt{-4\iota(\pm\infty)+\ve^2 H^2_\varSigma}\Big),
\end{equation}
and 
\[
\iota(\infty)=1-3(1+\sigma_\ve)^2, \quad \iota(-\infty)=1-3(-1+\sigma_\ve)^2.
\]
We can assume that the Wronskian at  $0$ is  $1$. If the right hand side of (\ref{match 5}) satisfies 
\begin{equation}
\int_\R g({\tt t}{, {\tt y})}  \partial_\st U(\st) e^{\,-\ve H_\varSigma
\st}\,d{\tt t}=0,\quad {\forall {\tt y}\in \varSigma} 
\label{ort 1}
\end{equation}
we can write $\varphi=\mathcal G(g)$ where
\begin{equation}
\label{lin 3}
\mathcal G(v)({\tt t}, {\tt y})=-\partial_\st U({\tt t})\int_0^{\tt t}
W(s)e^{\,-\ve H_\varSigma s} v(s,{{\tt y}})\,ds+W(\st)\int_{-\infty}^{\tt
t} \partial_s U(s)e^{\,-\ve H_\varSigma s} v(s, {{\tt y}})\,ds.
\end{equation}

Note that the orthogonality condition (\ref{ort 1}) guarantees that the function $\mathcal G(g)$ is exponentially decaying whenever $g$ is exponentially decaying (in $\st$).  To be more precise let us assume for instance that 
\[ 
|g(\sy, \st)|(\cosh \st)^\mu\leq C,
\]
with $\mu\in (\eta+\ve H_\varSigma, -\eta]$, where $\eta=\max\{\eta^+, \eta^-\}<0$. 
Then we have
\[
|\varphi(\sy, \st)|(\cosh \st)^\mu\leq C,
\]
as well.

Using this we can  determine the function $\psi_0$. To do this  we need to chose $h_0$ such that (\ref{ort 1}) is satisfied, in other words:
\[
\int_\R (\st+\ve h_0)\partial_\st U(\st)^2 e^{\,-\ve H_\varSigma \st}\,d\st=0,
\]
hence
\[
h_0=\frac{-\int_\R \st \partial_\st U(\st)^2 e^{\,-\ve H_\varSigma \st}\,d\st}{\ve \int_\R \partial_\st U(\st)^2 e^{\,-\ve H_\varSigma \st}\,d\st}=\mathcal O(1),
\]
With this choice we define:
\[
\psi_0(\sy, \st)=\mathcal G\big((\st+\ve h_0) \partial_\st U\big)|A_\varSigma|^2.
\]
Note that we have:
\[
(|\st|+\ve |h_0|)|\partial_\st U(\st)|(\cosh\st)^\mu< C,\quad 0<\mu<-\eta,
\]
and as a consequence 
\[
|\psi_0(\st,\sy)|(\cosh \st)^\mu<C,\quad 0<\mu<-\eta.
\]
Sometimes it is convenient to derive a  more refined estimate taking into account the fact that the leading order term on the right hand side is  $\st \partial_\st U(\st)=\mathcal O\big(|\st| (\cosh \st)^{\eta^\pm})\big)$, $|\st|\to \pm\infty$.  Thus, we consider (\ref{match 5}) assuming that 
\begin{equation}
\label{match 7}
|g(\sy, \st)|(1+|\st|)^\beta(\cosh \st)^\mu\leq C,
\end{equation}
with $\mu\in (\eta+\ve H_\varSigma, -\eta]$, and $\beta\in \R$. By a simple argument we have as well:
\begin{equation}\label{match 8}
|\mathcal G(g)(\st, \sy)|(1+|\st|)^\beta(\cosh \st)^\mu\leq C.
\end{equation}

\section{Delaunay solutions of the Cahn-Hilliard equation} 
\subsection{The Lyapunov-Schmidt reduction}\label{ls}
While our formal considerations in the above section were valid for any embedded CMC surface $\varSigma$ in $\R^d$ in what follows we will focus on a special example when $\varSigma=D_\tau$, i.e. it is a Delaunay unduloid. Since we are interested in functions which are periodic in the direction of the $x_d$ axis with the minimal period equal to that of $D_\tau$ we will introduce the manifold $\mathring  D_\tau$ which is obtained by identifying the set $D_\tau\cap \{x_d=0\}$ with the set $D_\tau\cap 
\{x_d=2T_\tau\}$. The set $\mathring D_\tau$ is homeomorphic to the $d-1$ dimensional torus $\mathbb T^{d-1}$.

Fist we note that the approximate solution $w=U+\ve^2\psi_0$ is so far only defined in $\mathcal N_\delta$, which is a tubular neighborhood of $\pD$. To extend $w$ to the whole space let us define
\begin{equation*}
 {\mathbb H}({\tt x})=\begin{cases}
           1+\sigma_\ve&\text{ if }{\tt x}\in\pD^+,\\
-1+\sigma_\ve&\text{ if }{\tt x}\in\pD^-,
          \end{cases}
\end{equation*}
where $\pD^+$, $\pD^-$ denote, respectively, the interior and the exterior of $\pD$.

Let us notice that the function $w$  approaches $\mathbb H$ exponentially. Indeed,  we have
\[
|Y^*_{\ve, h} w({\tt y}, {\tt t}) -Y_{\ve, h}^*{\mathbb  H}({\tt y}, {\tt
t})|\leq C_\mu e^{\,-\mu|{\tt t}|}, \quad \sy \in \pD, \quad{\tt t}\in\left
(-\frac{\delta}{\ve}, \frac{\delta}{\ve}\right),
\]
for any  $\mu\in (0, -\eta)$. To make the definition of $w^*$ precise we 
let $\chi$ to be a cutoff function such that $\chi(s)=1$, when $|s|\leq \frac{1}{2}$ and $\chi(s)=0$, when $|s|\geq 1$. Next, we define a cutoff  function $\chi^*$ supported in $\mathcal N_\delta$ by:
\[
Y^*_{\ve, h}\chi^*({\tt t})=\chi\left(\frac{{\ve\tt t}}{\delta}\right).
\]
We can define the global approximate solution $w^*$ by
\begin{equation}
\label{ls 1}
w^*({\tt x})= w({\tt x})\chi^*({\tt x})+{\mathbb H}({\tt x}) \big[1- \chi^*({\tt x})\big].
\end{equation}

Now we look for a solution of   the equation \eqref{CH} in the form 
\begin{equation*}
 u=w^*+\varphi
\end{equation*}
 where $\varphi$ is a small (in a way to be specified) function. 
 Thus our problem can be stated:
find  $\varphi\colon \R^{d-1}\times S^1_{2T_\tau}\to \R$, which is one-periodic
with period $2T_\tau$ in the $x_d$ direction,  such that
\begin{equation}
\label{ls 2}
N_\ve(w^*+\varphi)=\ell_\ve, \quad \mbox{in}\ \R^{d-1}\times S^1_{2T_\tau}, 
\end{equation}
where $S^{1}_{2T_\tau}$ is a circle  of radius $2T_\tau$,
\[
N_\ve(u)=\ve\Delta u+\frac{1}{\ve}u(1-u^2),
\]
and $\ell_\ve$ is the Lagrange multiplier defined in (\ref{eq U}). 
Let us recall that we want our solution to be rotationally symmetric. That is,  if by $\mathcal R_\theta$ we denote the rotation of $\R^d$ about the $x_d$ axis by angle $\theta$ then we should have:
\[
(w^*+\varphi)({\tt x})=(w^*+\varphi)(\mathcal R_\theta {\tt x}).
\]
Since we already have (by definition) 
\[
w^*({\tt x})=w^*(\mathcal R_\theta {\tt x}),
\]
than as a result we will have $\varphi({\tt x})=\varphi(\mathcal R_\theta {\tt x})$ as well, as can be seen easily from the proceeding construction. 

Since the function $\varphi$ appearing in  (\ref{ls 2}) is expected to be small it is natural to expand the nonlinear operator $N_\ve$ and write:
\[
L_\ve\varphi=-N_\ve(w^*)-Q_\ve(\varphi)+\ell_\ve
\]
where
\[
L_\ve\varphi=DN_\ve(w^*)\varphi, \quad Q_\ve(\varphi)=N(w^*+\varphi)-N_\ve(w^*)-L_\ve\varphi.
\]
The strategy, based on the Lyapunov-Schmidt reduction is clear. Indeed, we expect that due to the  $d$ dimensional bounded (and periodic) kernel of the Jacobi operator $\mathcal J_\tau$  which is associated  to translations of $D_\tau$ in the directions of the coordinate axis ${\tt e}_j$, $j=1, \dots, d$, the linear operator $L_\ve$ should have $d$ dimensional  kernel spanned, roughly speaking, by the functions $Z^{T,{\tt e}_j}_{\tau, \ve}$ where
\begin{equation}
\label{ls 3}
Y^*_{\ve, h}Z^{T, {\tt e}_j}_{\tau, \ve}({\tt y}, {\tt t})= {\tt V}(\sy, \st)\Phi_\tau^{T, {\tt e}_j}(\sy), \quad j=1, \dots, d, \quad {\tt V}(\sy, \st)=\partial_\st w=\partial_\st U(\st)+\ve^2\partial_\st \psi_0(\sy, \st).
\end{equation}
Notice also that  in general any function $Z_\ve$ such that
\[
Y^*_{\ve, h}Z_\ve({\tt y}, {\tt t})= \Phi({\tt y}) {\tt V}(\sy,{\tt t}), 
\]
is also ``almost" in the kernel of $L_\ve$, in the sense that $Y^*_{\ve, h}
(L_\ve Z_\ve)=o(1)$. We introduce a linear subspace of $L^2(\pD\times \R)$ of
functions that are orthogonal to $Z_\ve$ by:
\begin{equation}
\label{def XXX}
\mathcal X=\left\{\varphi\in L^2(\pD\times\R)\left| \int_\R \varphi(\sy, \st) {\tt V}(\sy, \st)\, d\st=0\right.\right\}.
\end{equation}
By $\Pi$ we denote  the orthogonal projection  on $\mathcal X$. We set $\varphi=\varphi^\parallel+\varphi^\perp$,  where
\[
Y^*_{\ve, h}\varphi^\parallel\in \mathcal X, \quad Y^*_{\ve, h}\varphi^\perp=Z{\tt V}\in \mathcal X^\perp.
\]
Finally we split our problem into two equations:
\begin{align}
\label{ls 4}
\Pi\circ Y^*_{\ve, h} [N_\ve(w^*+\varphi^\parallel+\varphi^\perp)-\ell_\ve]&=0, \\
\label{ls 5}
({\mathrm {Id}}\, -\Pi)\circ Y^*_{\ve, h}[N_\ve(w^*+\varphi^\parallel+\varphi^\perp)-\ell_\ve]&=0. 
\end{align}
When solving (\ref{ls 4}) we use the fact that the associated linear operator is
coercive on $\mathcal X$. To solve (\ref{ls 5})  we will make use of the  theory
of solvability of the Jacobi operator $\mathcal J_\pD$.  An additional, somewhat technical, step which we have omitted in this informal discussion   is to "transfer" the original problem from the space of functions defined on $ \R^{d-1}\times S^1_{2T_\tau}$ to a space of functions defined on $\pD\times \R$. We will explain these details  in  section \ref{sec glue} but first we will introduce and study a linear operator which is essentially the expression of  $L_\ve$ in the Fermi coordinates of $\pD$. 

%
%
\subsection{Linear theory for a model problem}\label{sec lin model}

In this section we will develop the necessary theory to deal with  the operator
$L_\ve$.  To this end we will consider the operator 
\begin{equation}
\label{lin 1}
\mathbb L_\ve=\ve \Delta_{\mathring  D_\tau}-(H_{\mathring  D_\tau}+\ve (\st+\ve h_0) \chi(\ve\st/\delta) |A_{\mathring D_\tau}|^2)\partial_\st+\ve^{-1} \partial_{\tt t}^2+\ve^{-1}f'(w),
\end{equation}
where $\chi(s)$ is a cutoff function supported in $(-1,1)$ and equal to $1$ in $(-1/2,1/2)$.  
Note that this operator is defined for functions $\phi\colon \mathring  D_\tau \times \R\to \R$ (and not just functions defined on $\mathring D_\tau\times (-\frac{\delta}{\ve}, \frac{\delta}{\ve})$).  It is clear that $Y^*_{\ve, h}L_\ve\approx \mathbb L_\ve$. Although the function $w=U+\ve^2\psi_0$ depends on both variables $(\sy, \st)$ in some sense the operator $\mathbb L_\ve$ separates variables. 
To see this, with  $\partial_\st w=\partial_\st (U+\ve^2\psi_0)={\tt V}$, we   consider functions of the form:
\[
\varphi(\sy, \st)={\tt V}(\sy, \st)Z(\sy).  
\]
Observe that, by  construction of $w=U+\ve^2\psi_0$, combining equations (\ref{eq U}) and (\ref{match 3}) multiplied by $\ve^2$ we get
\[
\ve^{-1} \partial_\st^2 w-(H_\pD +\ve (\st+\ve h_0) |A_{\mathring D_\tau}|^2)\partial_\st w+\ve^{-1} f(w)=\ell_\ve-\ve^3(\st+\ve h_0) |A_{\mathring D_\tau}|^2)\partial_\st\psi_0+\ve^{-1}Q(\ve\psi_0),
\]
where $Q(v)=f(U+v)-f(U)-f'(U) v$.
Differentiating  this equation in $\st$ we get  for ${\tt V}=\partial_\st  w$:
\begin{equation}
\begin{aligned}
\ve^{-1}\partial_\st^2 {\tt V}-(H_{\mathring D_\tau}+\ve(\st+\ve h_0) |A_{\mathring D_\tau}|^2)\partial_\st {\tt V}-\ve|A_{\mathring D_\tau}|^2 {\tt V}&+\ve^{-1} f'(w) {\tt V}=-\ve^3 |A_{\mathring D_\tau}|^2\partial_{\st}\psi_0\\
&\qquad-\ve^3(\st+\ve h_0) |A_{\mathring D_\tau}|^2\partial_\st^2\psi_0 -\ve^{-1}\partial_{\st} Q(\ve^2\psi_0).
\end{aligned}
\label{geom lin}
\end{equation}
From this, using the definition of $\mathbb L_\ve$ in (\ref{lin 1}) we get:
\begin{equation}\label{dschinn}
\mathbb L_\ve(\varphi)=\mathbb L_\ve({\tt V} Z)=\ve {\tt V}\mathcal J_{\mathring D_\tau}Z+B_\ve( Z)
\end{equation}
where $\mathcal J_{\mathring D_\tau}$ is the Jacobi operator on $\mathring D_\tau$ and
\begin{multline*}
B_\ve(Z)=-\ve({\tt t} +\ve h_0)(1-\chi^*)|A_\pD|^2 Z{\partial_{\tt t}}{\tt
V}+2\ve\nabla_\pD{\tt V}\cdot\nabla_\pD Z\\
+\ve Z\Delta_\pD{\tt V}-\ve^3\left[({\tt t}+\ve h_0)|A_\pD|^2\partial_{\tt
t}^2\psi_0+|A_\pD|^2\partial_{\tt t}\psi_0\right] Z\\
+\ve^{-1}\left[f'(w)\partial_\st w-f'(U)\partial_\st w
-\ve^2 f''(U)\partial_{\tt t} U\psi_0\right] Z.
\end{multline*}
We note that  
\begin{equation}
B_\ve(Z)=\mathcal O(\ve^3)\|Z\|_{\mathcal C^1(\pD)}.
\label{est b}
\end{equation}
Identity (\ref{geom lin}) and its consequence  (\ref{dschinn}) is the key calculation which allows to use the usual Lyapunov-Schmidt reduction scheme, as w explained in the introduction. Indeed, if we had taken as the approximate solution only the function $U$ then differentiating the equation (\ref{eq U}) for $U'=\partial_\st U$ we would have gotten
\[
\ve^{-1}\partial_\st^2 U'- H_\pD \partial_\st U'+\ve^{-1}f'(U) U'=0.
\]
This equation, unlike (\ref{geom lin}), does not carry any information about the geometry of $\pD$ besides its mean curvature which is constant. Following the method of \cite{MR2032110} or \cite{dkp_dg}, \cite{costa}, \cite{nested} we would have to perturb the surface $\pD$ additionally introducing new unknown functions in our problem. With the approach presented here this is no longer necessary and the  Lyapunov-Schmidt procedure in this version  is simpler.
Recalling that the linearization of the mean curvature operator is the Jacobi operator which depends on the second fundamental form, we see that  the operator $\mathbb L_\ve$ is naturally compatible with the geometric context of our problem.  To put it differently:  the operator $\mathbb L_\ve$ is, up to negligible terms,  the correct linearization of the Cahn-Hilliard operator near a solution whose zero level set is the constant curvature surface $\pD$.

%

To develop invertibility theory for $\mathbb L_\ve$ we will we employ two basic
facts. First, we observe that  on the  subspace:
\[
\mathcal Y:=\left\{\varphi(\sy, \st)={\tt V}(\sy, \st) Z(\sy)\mid\int_{\pD} Z(\sy) \Phi_\tau^{T, {\tt e}_j}(\sy)\, d\sy=0, j=1, \dots, d\right\},
\]
we have
\[
\langle {\mathbb L}_\ve\varphi, \varphi\rangle\geq C\ve \|\varphi\|^2_{L^2(\pD\times \R)}.
\] 
Second,  when we consider  $\varphi\in \mathcal X$ (space $\mathcal X$ is defined in (\ref{def XXX})) and $g\in L^2(\pD\times \R)$   such that $\varphi$ is a bounded   solution of the problem
\[
{\tt L}_\ve \varphi=g,\qquad {\tt L}_\ve =\ve\Delta_{\pD}+\ve^{-1}\partial^2_{\st\st}+\ve^{-1} f'(H),
\]
then we  have 
\begin{equation}
\label{key a priori}
\|\varphi\|_{L^2(\pD\times \R)}\leq C\ve\|g\|_{L^2(\pD\times \R)}.
\end{equation}
To prove this estimate we use a contradiction argument which  relies on  the fact that from 
\[
\int_\R |v'|^2-f'(H)v^2\geq C\|v\|^2_{L^2(\R)}, \quad \mbox{if}\
\int_\R v(t) H'(t)\,dt=0,
\]
it follows  that the bilinear form 
\[
{\tt B}_\ve (\varphi):=\langle {\tt L}_\ve\varphi, \varphi\rangle,
\]
is coercive on $\mathcal X$.   In the same way as (\ref{key a priori}) is shown one can  prove:
\begin{equation}
\label{key a priori 2}
\|\varphi\|_{H^1(\pD\times \R)}\leq C\|g\|_{L^2(\pD\times \R)}.
\end{equation}
The loss of the factor of $\ve$ in this last estimate is related with the scaling of the problem. We refer the reader to \cite{dkp_dg} or  \cite{nested} where results similar to estimates (\ref{key a priori})  and (\ref{key a priori 2}) were proven.

At the same time we can use (\ref{key a priori 2}) and the coercivity of the bilinear form ${\tt B}_\ve(\varphi)$ to solve the equation 
\begin{equation}
\label{sys lin 000}
\Pi_{\mathcal X}\mathbb L_\ve \varphi =g,
\end{equation}
where $\Pi_{\mathcal X}$ is the projection on $\mathcal X$, in $\mathcal X$. To do this we write $\mathbb L_\ve ={\tt L}_\ve+(\mathbb L_\ve-{\tt L}_\ve)$ and use a perturbation argument. The solution will still satisfy estimate (\ref{key a priori 2}).
{{The perturbation argument is as follows: for $g\in \mathcal X$ we solve
\[
{\tt L}_\ve\phi=g+cH',
\]
where $c=\frac{\int g H'\,dt}{\int (H')^2\,dt}$. Then we define a map
\[
G_{\mathcal X}(g)=\phi-V\frac{\int \phi V\,dt}{\int V^2\,dt}.
\]
Note that $G_\mathcal X\colon \mathcal X\to \mathcal X$. Next we check
\[
\|\Pi_{\mathcal X} \mathbb L_\ve G_{\mathcal X}(g)-g\|_{L^2(\pD\times \R)}\leq o(1)\|g\|_{L^2(\pD\times \R)}.
\]

Therefore $\Pi_{\mathcal X} \mathbb L_\ve G_{\mathcal X}$ is invertible as a map from $\mathcal X$ to itself and   we can  define
\[
(\Pi_{\mathcal X}\mathbb L_\ve)^{-1}=G_{\mathcal X}(\mathbb L_\ve G_{\mathcal X})^{-1}.
\]
Moreover it is rather straightforward to show that a solution to (\ref{sys lin 000}) will satisfy estimate (\ref{key a priori 2}). 
}}

We will use  these observations  to solve the following model equation:
\begin{equation}
\label{lin 2}
\mathbb L_\ve \varphi= g(\sy, \st),
\end{equation}
where we will assume initially that $g\in L^p(\pD\times \R)$. We look for a solution in the form $\varphi=\varphi^\parallel+\varphi^\perp$,  where
\[
\varphi^\parallel\in \mathcal X, \quad \varphi^\perp=Z{\tt V}\in \mathcal X^\perp.
\]
{{We write
\[
\mathbb L_\ve \varphi=\Pi_{\mathcal X} \mathbb L_\ve (\varphi^\parallel+\varphi^\perp)+\Pi_{\mathcal X^\perp}\mathbb L_\ve (\varphi^\parallel+\varphi^\perp),
\]
and then we need to solve
\[
\begin{aligned}
\Pi_{\mathcal X} \mathbb L_\ve \varphi^\parallel  &= \Pi_{\mathcal X} g-\Pi_{\mathcal X} \mathbb L_\ve \varphi^\perp\\
\Pi_{\mathcal X^\perp}\mathbb L_\ve \varphi^\perp&=\Pi_{\mathcal X^\perp} g- \Pi_{\mathcal X^\perp}\mathbb L_\ve \varphi^\parallel.
\end{aligned}
\]
The idea is that terms $\Pi_{\mathcal X} \mathbb L_\ve \varphi^\perp$ and $\Pi_{\mathcal X^\perp}\mathbb L_\ve \varphi^\parallel$ are of smaller order because $\mathbb L_\ve V=o(1)$ so that the coupling between the two equations is rather weak. 
Another important point is that 
\[
\mathbb L_\ve \varphi^\perp=\ve V{\mathcal J}_\pD Z + B_\ve(Z),
\]
where $B_\ve(Z)$ is small (see (\ref{est b})). 
We decompose accordingly $g=g^\parallel+g^\perp$ , $g^\perp=\Xi{\tt V}$ and  $B_\ve(Z)=B_\ve^\parallel(Z)+ B_\ve^\perp(Z)$, $B^\perp_\ve(Z)=\Upsilon_\ve(Z) {\tt V}$ and look for a solution of the system:
\begin{equation}
\begin{aligned}
\Pi_{\mathcal X}\mathbb L_\ve\varphi^\parallel&=g^\parallel -B_\ve^\parallel (Z), \\
\ve {\mathcal J}_\pD Z+\Upsilon_\ve(Z)&=\Xi- \int_\R(\mathbb L_\ve \varphi^\parallel){\tt V}\,d\st+\sum_{j=1}^{d} c_j \Phi^{T, {\tt
e}_j}_{\tau}.
\end{aligned}
\label{sys lin 1}
\end{equation}
Note that in the second equation we have introduced Lagrange multipliers  ${c}_j$ to be determined. 

Although the two equations in (\ref{sys lin 1}) are coupled but this coupling is weak and we can solve the system without any difficulty using invertibility of $\Pi_{\mathcal X}$ and ${\mathcal J}_\pD$,  a fixed point argument and estimates (\ref{est b}),  (\ref{key a priori}) and (\ref{key a priori 2}). We leave this to the reader.  

\medskip

Given that we can solve  (\ref{sys lin 1})  our purpose is to find suitable estimates for the solution of  the problem  
\[
\mathbb L_\ve\varphi =g(\sy, \st)
\]
on $\mathcal X$  assuming that  
\[
g(\sy, \st)=\mathcal O(e^{\,-\mu|\st|}), \quad |\st|\to \infty.
\]
In particular we would like to know that $\varphi(\st, \sy)=\mathcal O(e^{\,-\mu|\st|})$ as well. 
This is straightforward by comparison principle once we know for example that $\varphi$ is bounded. Thus the main issue is to obtain $L^\infty$ control for $\varphi$.  We will go a little further now and show how to control  {\it a priori} certain weighted H\"older norms of $\varphi^\parallel$ and $\varphi^\perp$ (see decomposition in (\ref{sys lin 1})).  
%
%

In order to simplify the argument and avoid keeping track of negative powers of $\ve$ appearing on the right hand side of various estimates we will rescale the $\sy$ variable. Thus we will introduce:
\[
\tilde \sy =\frac{\sy}{\ve}, \quad \tilde \st=\st.
\]
We will denote $\mathring D_{\tau, \ve}=\frac{1}{\ve}\pD$.
We consider the manifold $\mathring D_{\tau, \ve}\times \R$ equipped with the product metric and  the associated Levi-Civita connection. We define weighted H\"older norms on $\mathring D_{\tau, \ve}\times \R$:
\begin{equation}\label{holder 1}
\begin{aligned}
\|u\|_{\mathcal C^{0,\alpha}_\mu(\mathring D_{\tau, \ve}\times
\R)}&=\sup_{\tilde \st\in \R}(\cosh \tilde \st)^\mu \|u\|_{\mathcal
C^{0,\alpha}(\mathring D_{\tau, \ve}\times (\tilde\st -1, \tilde\st+1))},\\
{\|u\|_{\mathcal C^{1,\alpha}_\mu(\mathring D_{\tau,\ve}\times\R)}}
&={\|u\|_{\mathcal C^{0,\alpha}_\mu(\mathring D_{\tau,\ve}\times
\R)}+\|\nabla_{\mathring D_{\tau,\ve}\times \R}u\|_{\mathcal
C^{0,\alpha}_\mu(\mathring D_{\tau,\ep}\times \R)},}\\
\|u\|_{\mathcal C^{2,\alpha}_\mu(\mathring D_{\tau, \ve}\times
\R)}&=\|u\|_{\mathcal C^{0,\alpha}_\mu(\mathring D_{\tau, \ve}\times
\R)}+\|\nabla_{\mathring D_{\tau, \ve}\times \R}u\|_{\mathcal
C^{0,\alpha}_\mu(\mathring D_{\tau, \ve}\times \R)}+\|\nabla^2_{\mathring
D_{\tau, \ve}\times \R}u\|_{\mathcal C^{0,\alpha}_\mu(\mathring D_{\tau,
\ve}\times \R)}.
\end{aligned}
\end{equation}
Above  $\nabla_{\mathring D_{\tau, \ve}\times \R}$ and $\nabla^2_{\mathring
D_{\tau, \ve}\times \R}$ respectively denote the gradient and the Hessian on the
manifold $\mathring D_{\tau, \ve}\times \R$. 


Given functions  $u$ and $g$ on $\colon \pD\times \R$ and  we introduce functions 
\[
\tilde u(\tilde\sy, \tilde\st)=u(\ve\tilde\sy, \tilde \st), \quad \tilde w(\tilde\sy, \tilde\st)=w(\ve\tilde\sy, \tilde \st),\quad \tilde g(\tilde\sy, \tilde\st)=\ve g(\ve\tilde \sy, \tilde\st).
\]
and set
\[
\tilde {\mathbb L}_\ve\tilde u=\Delta_{\mathring D_{\tau, \ve}} \tilde u+\partial^2_{\tilde\st}\tilde u+f'(\tilde w)\tilde u+\ve\tilde q\partial_{\tilde\st}\tilde u,
\]
where 
\[
\tilde q := H_{\mathring  D_\tau}+\ve (\tilde \st+\ve h_0) \chi(\ve\tilde \st/\delta) |A_{\mathring D_\tau}(\ve \tilde \sy)|^2
\] 
is a bounded function. The linear problem we  consider  is:
\begin{equation}
\label{lin 5.1.2}
\Pi_{\mathcal X}\tilde{\mathbb L}_\ve\tilde u=\tilde g, \quad \mbox{in}\  \mathring D_{\tau, \ve}\times \R,
\end{equation}
where now we assume simply that the right hand side satisfies the orthogonality condition (c.f. (\ref{lin 3})):
\[
\int_{\R}\tilde g(\tilde \sy, \tilde\st) {\tt V}(\ve\tilde \sy, \tilde\st)\, d\tilde\st=0.
\]
For this problem we can derive {\it a priori} estimates using the method of   \cite{dkp_dg}, \cite{costa}, \cite{nested}.  From this, by a perturbation  argument, we will be able get corresponding estimates for $\varphi^\parallel$ satisfying the first equation in (\ref{sys lin 1}). Let us explain briefly the main steps.

\medskip
\noindent
%
%
\medskip
\noindent
{\bf Step 1.}
We consider a problem of the form:
\begin{equation}
\Delta_{y} \phi+\partial^2_t\phi+f'(H(t))\phi=0, \quad \mbox{in}\ \R^{d-1}\times \R.
\label{lin 5.2}
\end{equation}
The following result is known (\cite{dkp_dg}, \cite{costa}):
\begin{lemma}\label{nondeg H}
Let $\phi$ be a bounded solution of (\ref{lin 5.2}). Then $\phi=c H'(\st)$, with some constant $c$.
\end{lemma}

\medskip
\noindent
{\bf Step 2.} 
Consider now equation  (\ref{lin 5.1.2}) and assume that $\tilde g\in \mathcal C^{0,\alpha}_\mu (\mathring D_{\tau, \ve}\times \R)$, with $\mu\in (0, |\eta|)$. The we have:
\begin{lemma}\label{step 2}
There exists a constant $C>0$ such that for all sufficiently small $\ve$ any bounded  solution of (\ref{lin 5.1.2}) satisfies:
\begin{equation}
\label{lin 5.3}
\|\tilde u\|_{\mathcal C^{2,\alpha}_\mu(\mathring D_{\tau, \ve}\times \R)}\leq C \|\tilde g\|_{\mathcal C^{0,\alpha}_\mu(\mathring D_{\tau, \ve}\times \R)}.
\end{equation}
\end{lemma}
A proof of this lemma, which relies on Lemma \ref{nondeg H} and a contradiction argument,  follows the same lines as the proof of Lemma  5.2 in  \cite{nested} (see also similar results in \cite{dkp_dg}, \cite{costa}).

Now we should go back to the original variables. We define weighted H\"older norms on $\pD\times \R$ similarly as in (\ref{holder 1}):
\begin{equation}\label{holder 2}
\begin{aligned}
\|u\|_{\mathcal C^{0,\alpha}_\mu(\pD\times \R)}&=\sup_{\st\in \R}(\cosh \st)^\mu
\|u\|_{\mathcal C^{0,\alpha}(\mathring D_{\tau}\times (\st -1, \st+1))},\\
{\|u\|_{\mathcal C^{1,\alpha}_\mu(\mathring D_{\tau}\times
\R)}}&={\|u\|_{\mathcal
C^{0,\alpha}_\mu(\mathring D_{\tau}\times \R)}+\|\nabla_{\mathring
D_{\tau}\times \R}u\|_{\mathcal C^{0,\alpha}_\mu(\mathring D_{\tau}\times
\R)}}\\
\|u\|_{\mathcal C^{2,\alpha}_\mu(\mathring D_{\tau}\times \R)}&=\|u\|_{\mathcal
C^{0,\alpha}_\mu(\mathring D_{\tau}\times \R)}+\|\nabla_{\mathring
D_{\tau}\times \R}u\|_{\mathcal C^{0,\alpha}_\mu(\mathring D_{\tau}\times
\R)}+\|{\nabla}^2_{\mathring D_{\tau}\times \R}u\|_{\mathcal
C^{0,\alpha}_\mu(\mathring D_{\tau}\times \R)}.
\end{aligned}
\end{equation}
We note that if for a given function $u\colon \mathring D_\tau\times \R\to \R$ we set $\tilde u(\tilde \sy, \tilde \st)=u(\ve\tilde \sy, \tilde \st)$ then
{{we have
\begin{equation}
\label{relation norms}
 \|\tilde u\|_{\mathcal C^{\ell,\alpha}_\mu(\mathring D_{\tau_\ve}\times \R)}=
\sum_{0\leq k+m\leq \ell} \ve^{m}\|\partial^k_{\st} D^m_\pD u\|_{\mathcal
C^{0}_\mu(\mathring D_{\tau}\times \R)}+\sum_{0\leq k+m\leq \ell}
\ve^{m+\alpha}[\partial^k_{\st} D^m_\pD u]_{\alpha, \mu, \mathring
D_{\tau}\times \R},
\end{equation}
where $[\cdot]_{\alpha, \mu, \mathring D_{\tau}\times \R}$ is the weighted H\"older  seminorm.
Consequently by $\mathcal E^{\ell,\alpha}_\mu(\mathring D_{\tau}\times \R)$ we
denote the space of functions on $\pD\times \R$  with the norm
\[
\|u\|_{\mathcal E^{\ell,\alpha}_\mu(\mathring D_{\tau}\times \R)}:=
\sum_{0\leq k+m\leq \ell} \ve^{m}\|\partial^k_{\st} D^m_\pD u\|_{\mathcal
C^{0,\alpha}_\mu(\mathring D_{\tau}\times \R)}.
\]
With this definition we have 
\[
\|u\|_{\mathcal E^{0,\alpha}_\mu(\mathring D_{\tau}\times \R)}=\|u\|_{\mathcal C^{0,\alpha}_\mu(\mathring D_{\tau}\times \R)},
\]
while with the notation of (\ref{relation norms}) we get
\[
C^{-1} \|\tilde u\|_{\mathcal C^{\ell,\alpha}_\mu(\mathring D_{\tau_\ve}\times \R)}\leq \|u\|_{\mathcal E^{\ell,\alpha}_\mu(\mathring D_{\tau}\times \R)}\leq  C \ve^{-\alpha}\|\tilde u\|_{\mathcal C^{\ell,\alpha}_\mu(\mathring D_{\tau_\ve}\times \R)}.
\]
}} 
From this and Lemma \ref{step 2} it follows for $\ell =0, 1, 2$:
\begin{equation}\label{holder 3}
\|u\|_{\mathcal E^{\ell,\alpha}_\mu(\mathring D_{\tau}\times \R)}\leq C\ve^{1-\alpha}\|g\|_{\mathcal E^{0,\alpha}_\mu(\mathring D_{\tau}\times \R)}.
\end{equation}
 The procedure described above shows that  we can  control the size of $\mathcal E^{\ell,\alpha}_\mu(\mathring D_{\tau}\times \R)$ norm of $\varphi^\parallel$ in (\ref{sys lin 1}) obtaining:
 \begin{equation}\label{holder 4}
 \|\varphi^\parallel\|_{\mathcal E^{\ell,\alpha}_\mu(\mathring D_{\tau}\times \R)}\leq C\ve^{1-\alpha}(\|g^\parallel\|_{\mathcal E^{0,\alpha}_\mu(\mathring D_{\tau}\times \R)}+\ve^2\|g^\perp\|_{\mathcal E^{0,\alpha}_\mu(\mathring D_{\tau}\times \R)}).
 \end{equation}

Finally, we notice that for the second equation in (\ref{sys lin 1}) using elliptic theory we can get   H\"older estimates and since $\varphi^\perp=Z{\tt V}$ we find:
\begin{equation}\label{holder 5}
 \|\varphi^\perp\|_{\mathcal C^{\ell,\alpha}_\mu(\mathring D_{\tau}\times \R)}\leq C\ve^{-1}\|g^\perp\|_{\mathcal C^{0,\alpha}_\mu(\mathring D_{\tau}\times \R)}.
 \end{equation}
 
 }}

\subsection{The linear problem in the whole space}\label{sec glue}
Now we will use the theory outlined above to solve the following problem:
\begin{equation}
\label{lin 10}
\ve \Delta \varphi +\frac{1}{\ve} f'(w^*)\varphi=g({\tt x}), \quad \mbox{in}\  \R^{d-1}\times S_{2T_\tau},
\end{equation}
From what we have said above it is in general not possible to find a solution with a reasonably bounded norm unless the right hand side satisfies some extra conditions, or equivalently, we need to introduce   Lagrange multipliers that correspond to natural invariances of the problem. Thus, we will solve
\begin{equation}
\label{lin 11}
\ve \Delta \varphi +\frac{1}{\ve} f'(w^*)\varphi=g({\tt
x})+\chi^*\sum_{j=1}^d{\tt c}_j Z^{T, {\tt e}_j}_{\tau, \ve}, \quad \mbox{in}\ 
\R^{d-1}\times S_{2T_\tau}
\end{equation}
where 
\[
Y_{\ve, h}^* \chi^*(\st)=\chi(\ve\st/\delta),   
\quad Y^*_{\ve, h}Z^{T, {\tt e}_j}_{\tau, \ve}({\tt y}, {\tt t})= {\tt V}(\sy, \st)\Phi_\tau^{T, {\tt e}_j}(\sy), \quad j=1, \dots, d
\]

The idea is to solve (\ref{lin 11}) by {\it gluing} a solution defined near $\pD$ and another one defined away from $\pD$. To describe this construction rigorously we need some preparation. 
We introduce the function $q({\tt x})$ as follows:
\[
q({\tt x})=\begin{cases}f'(1+\sigma_\ve), \quad \mathrm{dist}\, ({\tt x}, D_\tau)>\delta/2,\\
f'(-1+\sigma_\ve), \quad \mathrm{dist}\, ({\tt x}, D_\tau)<-\delta/2,
\end{cases}
\]
and otherwise $q({\tt x})$ is a smooth function such that $\min\{f'(1+\sigma_\ve), f'(-1+\sigma_\ve)\}<q({\tt x})\leq \max\{f'(1+\sigma_\ve), f'(-1+\sigma_\ve)\}$. Note that $q({\tt x})=-2+\mathcal O(\ve)$.   
Finally, we need another cutoff function $\widetilde \chi$ such that $\widetilde\chi\chi^*=\chi^*$ (take for instance $Y_{\ve, h}^* \widetilde\chi(\st)=\chi(\ve\st/2\delta)$ and chose $\delta$ smaller so that the Fermi coordinates are defined in $\mathcal N_{2\delta}$). We want to find a solution of (\ref{lin 11}) in the form $\varphi=\chi^*\check\varphi\circ Y_{\ve, h} +\psi$, where the function $\check\varphi$ solves:
\begin{equation}
\label{lin 12}
{\mathbb L}_\ve \check \varphi=\Big(\widetilde\chi\{g+\chi^*\sum_{j=1}^d{\tt c}_j Z^{T, {\tt e}_j}_{\tau,
\ve}+({\mathbb L}_\ve-L_\ve)\check \varphi-[\chi^*, L_\ve]\check
\varphi+\ve^{-1}[q-f'(w^*)]\psi\}\Big),\quad \mbox{in} \
\pD\times \R,
\end{equation}
and the function $\psi$ solves
\begin{equation}
\label{lin 13}
{\ve\Delta\psi+
\ve^{-1}[(1-\chi^*)f'(w^*)+\chi^*q]\psi=(1-\chi^*)\{g+\chi^*\sum_{j=1}^d{\tt
c}_jZ_{\tau,\ve}^{T,{\tt e}_j}\\-[\chi^*, L_\ve]\check\varphi\}
}
, \quad \mbox{in}\ \R^{d-1}\times S_{2T_\tau}.
\end{equation}
It  is clear that multiplying (\ref{lin 12}) by $\chi^*$ and adding the  equations (\ref{lin 12})--(\ref{lin 13}) and using the fact that $\widetilde\chi\chi^*=\chi^*$ we get the solution to our problem.  In the above and in what follows we abuse slightly notation writing for instance $\check \varphi$ as a function defined on $\pD\times \R$ and as a function defined on $\R^{d-1}\times S_{2T_\tau}$. It is understood that in the latter case we take $\check\varphi\circ Y_{\ve, h}$. To avoid complicated notions we will omit the composition with $Y_{\ve, h}$ or $Y^{-1}_{\ve, h}$ whenever it does not cause confusion. Thus the commutator $[\chi^*, L_\ve]\check \varphi$ in $\R^{d-1}\times
S_{2T_\tau}$ is 
\[ 
[\chi^*, L_\ve]\check\varphi=2\ve\nabla\check\varphi\circ Y_{\ve, h}\nabla\chi^*+\ve\check\varphi\circ Y_{\ve, h}\Delta\chi^*,
\] 
while in $\pD\times R$ we have to first express $L_\ve$ in local coordinate $(\sy, \st)$ (written as $Y^*_{\ve, h}L_\ve$) and calculate   $[\chi^*, Y^*_{\ve, h}L_\ve]\check\varphi$.

The function $g$ on the right hand side of this equation satisfies the following general assumptions on its asymptotic behaviour:{{
\begin{equation}\label{g hype}
\begin{aligned}
\|(g\chi^*)^\parallel\|_{\mathcal E^{0,\alpha}_\mu (\mathring D_{\tau}\times \R)}&\leq C,\\
\|(g\chi^*)^\perp\|_{\mathcal C^{0,\alpha}_\mu (\mathring D_{\tau}\times \R)}&\leq C\\
\|(1-\chi^*) g\|_{\mathcal C^{0,\alpha}(\R^{d-1}\times S_{2T_\tau})}&\leq C e^{\,-c_0/\ve}.
\end{aligned}
\end{equation}
}}
In addition we assume that $g$ is rotationally symmetric about the $x_d$ axis, namely if by $\mathcal R_\theta$ we denote the rotation of $\R^d$ about the $x_d$ axis by angle $\theta$ then $g(\mathcal R_\theta {\tt x})=g({\tt x})$.

In order to solve this coupled system we need to make sure that all terms on the right hand side that involve $\check\varphi$ and $\psi$ are small in suitable weighted H\"older and H\"older norms respectively. It is at this point that we need to chose the parameter $\delta$ in the definition of the tubular neighbourhood $\mathcal N_\delta$ small and dependent on $\ve$. Thus we take $\delta(\ve)=\ve^{2/3}$. This means in particular that for ${\tt x}\in \mathcal N_\delta$ we have $\ve \st({\tt x})=\mathcal O(\ve^{2/3})$. For reasons that will become clear soon we will also chose the H\"older exponent $\alpha$ in the definition of  $\mathcal C^{0,\alpha}_\mu (\mathring D_{\tau}\times \R)$ and $\mathcal C^{0,\alpha} (\mathring D_{\tau}\times \R)$ to be in the interval $(0,\frac{1}{10})$.  Finally, the parameter $\mu$ will be always taken in the interval $(0, |\eta|)$. 

Considering equation (\ref{lin 13}) we have the following:
\begin{lemma}
If $u$ is a solution of
\[{
\ve\Delta u-\frac1\ve [(1-\chi^*)f'(w^*)+\chi^*q] u=g}, \quad \mbox{in}\
\R^{d-1}\times S_{2T_\tau}, 
\]
then we have an {\it a priori} estimate:
\begin{equation}\label{lin 13 a}
\|u\|_{\mathcal C^{\ell, \alpha}(\R^{d-1}\times S_{2T_\tau})}\leq C\ve^{1-\ell-\alpha}\|g\|_{\mathcal C^{0, \alpha}(\R^{d-1}\times S_{2T_\tau})}.
\end{equation}
\end{lemma}
The proof of this lemma is  straightforward and it is omitted, for similar
results see for instance \cite{MR2646302}.
From this we get readily an {\it a priori} estimate for (\ref{lin 13}):
\begin{equation}\label{lin 13 b}
{\|\psi\| _{\mathcal C^{\ell, \alpha}(\R^{d-1}\times S_{2T_\tau})}\leq
C\ve^{1-\ell-\alpha}\|(1-\chi^*)g\|_{\mathcal C^{0, \alpha}(\R^{d-1}\times
S_{2T_\tau})}+e^{\,-c\ve^{-1/4}}(\sum_{j=1}^d|{\tt
c}_j|+\|\check\varphi\|_{\mathcal C^{1, \alpha}_\mu(\pD\times
\R)}).}
\end{equation}
From the theory developed in the previous section we can also obtain an {\it a priori} estimate for (\ref{lin 12}). If we write
\[
\mathfrak g=\widetilde\chi\{g+\chi^*\sum_{j=1}^d{\tt c}_j Z^{T, {\tt e}_j}_{\tau, \ve}+({\mathbb L}_\ve-L_\ve)\check \varphi-[\chi^*, L_\ve]\check \varphi+\ve^{-1}[q-f'(w^*)]\psi\},
\]  
then we have using (\ref{holder 4})
\begin{equation}
\label{parallel estimate}
\|\check \varphi^\parallel\|_{\mathcal E^{\ell, \alpha}_\mu(\pD\times \R)}\leq
C\ve^{{1}-\alpha}(\|\mathfrak g^\parallel\|_{\mathcal E^{0,
\alpha}_\mu(\pD\times \R)}+\ve^2\|\mathfrak g^\perp\|_{\mathcal C^{0,
\alpha}_\mu(\pD\times \R)}),
\end{equation}
and using (\ref{holder 5}):
\begin{equation}\label{perp estimate}
\|\check \varphi^\perp\|_{\mathcal C^{\ell, \alpha}_\mu(\pD\times \R)}\leq C\ve^{-1}\|\mathfrak g^\perp\|_{\mathcal C^{0, \alpha}_\mu(\pD\times \R)}.
\end{equation}
{{
We note that the weighted norms we use for $\check\phi^\parallel$ and $\check\phi^\perp$ are scaled differently with $\ve$. This slight nuisance is a result of our choice of the original scaling of the Cahn-Hilliard equation. We observe as well that with our definitions  $\|\cdot\|_{\mathcal E^{0, \alpha}_\mu(\pD\times \R)}= \|\cdot\|_{\mathcal C^{0, \alpha}_\mu(\pD\times \R)}$. 

We will now estimate $\|\mathfrak g^\parallel\|_{\mathcal E^{0, \alpha}_\mu(\pD\times \R)}$. To do this we observe that, with $\check \phi=\check\phi^\parallel+\check\phi^\perp$, we have:
\[
\begin{aligned}
\big\|\big(\widetilde\chi ({\mathbb L}_\ve-L_\ve)\check
\varphi\big)^\parallel\big\|_{\mathcal C^{0, \alpha}_\mu(\pD\times \R)}&\leq
C\delta(\ve)(\ve^{-1}\|\check \varphi^\parallel\|_{\mathcal E^{2,
\alpha}_\mu(\pD\times \R)}+ \ve\|\check \varphi^\perp\|_{\mathcal C^{2,
\alpha}_\mu(\pD\times \R)})\\
\|\big(\widetilde\chi [\chi^*, L_\ve]\check
\varphi\big)^\parallel\|_{\mathcal E^{0, \alpha}_\mu(\pD\times \R)}&\leq 
C\delta(\ve)^{-1}(\|\check \varphi^\parallel\|_{\mathcal E^{1,
\alpha}_\mu(\pD\times \R)}+\mathcal O(e^{\,-c\ve^{-1/3}})\|\check
\varphi^\perp\|_{\mathcal C^{1, \alpha}_\mu(\pD\times \R)}),\\
\big\|\big(\widetilde\chi\ve^{-1}[q-f'(w^*)]\psi\big)^\parallel\big\|_{\mathcal E^{0, \alpha}_\mu(\pD\times \R)}&\leq C\ve^{-1}\|\psi\|_{\mathcal C^{0, \alpha}(\R^{d-1}\times S_{2T_\tau})}.
\end{aligned}
\]
Next we estimate the orthogonal complement of these functions
\[
\begin{aligned}
\big\|\big(\widetilde\chi ({\mathbb L}_\ve-L_\ve)\check \varphi\big)^\perp\big\|_{\mathcal E^{0, \alpha}_\mu(\pD\times \R)}&\leq C(\|\check \varphi^\parallel\|_{\mathcal E^{2, \alpha}_\mu(\pD\times \R)}+ \ve^2\|\check \varphi^\perp\|_{\mathcal C^{2, \alpha}_\mu(\pD\times \R)})\\
\|\big(\widetilde\chi [\chi^*, L_\ve]\check \varphi\big)^\perp\|_{\mathcal E^{0, \alpha}_\mu(\pD\times \R)}&\leq  \mathcal O(e^{\,-c\ve^{-1/3}})\big(\|\check \varphi^\parallel\|_{\mathcal E^{1, \alpha}_\mu(\pD\times \R)}+\|\check \varphi^\perp\|_{\mathcal C^{1, \alpha}_\mu(\pD\times \R)}\big),\\
\big\|\big(\widetilde\chi\ve^{-1}[q-f'(w^*)]\psi\big)^\perp\big\|_{\mathcal C^{0, \alpha}_\mu(\pD\times \R)}&\leq C\ve^{-1}\|\psi\|_{\mathcal C^{0, \alpha}(\R^{d-1}\times S_{2T_\tau})}
\end{aligned}
\]
We can estimate the parameters ${\tt c}_j$ by projection of $\mathfrak g$ onto $Z^{T, {\tt e}_j}_{\tau, \ve}$.  Using the above estimate we get: 
\[
\begin{aligned}
|{\tt c}_j|&\leq C\big\{\|(\widetilde\chi\ g)^\perp\|_{\mathcal C^{0,
\alpha}_\mu(\pD\times \R)}
+\big\|\big(\widetilde\chi ({\mathbb L}_\ve-L_\ve)\check \varphi\big)^\perp\big\|_{\mathcal E^{0, \alpha}_\mu(\pD\times \R)}
\\
&\qquad +\|\big(\widetilde\chi [\chi^*, L_\ve]\check \varphi\big)^\perp\|_{\mathcal E^{0, \alpha}_\mu(\pD\times \R)}
+\big\|\big(\widetilde\chi\ve^{-1}[q-f'(w^*)]\psi\big)^\perp\big\|_{\mathcal C^{0, \alpha}_\mu(\pD\times \R)}\big\}
\\
&\leq C\left\{\|(\widetilde\chi\ g)^\perp\|_{\mathcal C^{0,
\alpha}_\mu(\pD\times \R)}
+\|\check \varphi^\parallel\|_{\mathcal E^{2, \alpha}_\mu(\pD\times \R)}+ \ve^2\|\check \varphi^\perp\|_{\mathcal C^{2, \alpha}_\mu(\pD\times \R)}
\right.\\
&\qquad 
\left.+
O(e^{\,-c\ve^{-1/3}})\big(\|\check \varphi^\parallel\|_{\mathcal E^{1, \alpha}_\mu(\pD\times \R)}+\|\check \varphi^\perp\|_{\mathcal C^{1, \alpha}_\mu(\pD\times \R)}\big)
+\ve^{-1}\|\psi\|_{\mathcal C^{0, \alpha}(\R^{d-1}\times S_{2T_\tau})}
\right\}
\\
&\leq C\left\{\|(\widetilde\chi\ g)^\perp\|_{\mathcal C^{0,
\alpha}_\mu(\pD\times \R)}+\|\check \varphi^\parallel\|_{\mathcal E^{2, \alpha}_\mu(\pD\times \R)}+ \ve^2\|\check \varphi^\perp\|_{\mathcal C^{2, \alpha}_\mu(\pD\times \R)}+\ve^{-1}\|\psi\|_{\mathcal C^{0, \alpha}(\R^{d-1}\times S_{2T_\tau})}\right\}
\end{aligned}
\]
Now we use estimates (\ref{parallel estimate})--(\ref{perp estimate}). After rearranging terms suitably and  using $\ve^{1-\alpha}\delta^{-1}(\ve)=o(1)$ to absorb $\check\varphi^\parallel$ in the first inequality below  we get 
\begin{equation}\label{a priori 1}
\begin{aligned}
\|\check \varphi^\parallel\|_{\mathcal E^{2, \alpha}_\mu(\pD\times \R)}&\leq
C\ve^{1-\alpha}\big\{\|(\widetilde\chi\ g)^\parallel\|_{\mathcal
E^{0, \alpha}_\mu(\pD\times \R)}  + \ve^2\|(\widetilde\chi\
g)^\perp\|_{\mathcal C^{0, \alpha}_\mu(\pD\times \R)}\\
&\qquad\qquad +\ve^{-1}\|\psi\|_{\mathcal
C^{0, \alpha}(\R^{d-1}\times S_{2T_\tau})}+{\delta(\ve) \|\check
\varphi^\perp\|_{\mathcal C^{2, \alpha}_\mu(\pD\times \R)}}\big\},
\\
\|\check \varphi^\perp\|_{\mathcal C^{2, \alpha}_\mu(\pD\times \R)}&\leq
C\ve^{-1}\big\{\|(\widetilde\chi\ g)^\perp\|_{\mathcal C^{0,
\alpha}_\mu(\pD\times \R)}+{
\|\check \varphi^\parallel\|_{\mathcal E^{2, \alpha}_\mu(\pD\times \R)}}
+\ve^{-1}\|\psi\|_{\mathcal C^{0,
\alpha}(\R^{d-1}\times S_{2T_\tau})}\big\}.
\end{aligned}
\end{equation}
From (\ref{lin 13 b}) we get as well for $\ell=0,1,2$:
\begin{equation}\label{lin 14}
\begin{aligned}
\|\psi\| _{\mathcal C^{\ell, \alpha}(\R^{d-1}\times S_{2T_\tau})}& \leq C\ve^{1-\ell-\alpha}\{\|(1-\chi^*)g\|_{\mathcal C^{0, \alpha}(\R^3)}\\
&\quad +e^{\,-c\ve^{-1/8}}(\|\widetilde \chi g\|_{\mathcal C^{0, \alpha}_\mu(\pD\times \R)}+\|\check\varphi^\parallel\|_{\mathcal E^{2, \alpha}_\mu(\pD\times \R)}+\|\check\varphi^\perp\|_{\mathcal C^{2, \alpha}_\mu(\pD\times \R)}+\|\psi\|_{\mathcal C^{0, \alpha}(\R^{d-1}\times S_{2T_\tau})})\}.
\end{aligned}
\end{equation}
Using the fact that $\delta(\ve)\ve^{-\alpha}=o(1)$  to absorb term  ${\delta(\ve) \|\check
\varphi^\perp\|_{\mathcal C^{2, \alpha}_\mu(\pD\times \R)}}$ appearing on the right hand side of the first inequality in (\ref{a priori 1}) and combining these estimates we get
\begin{equation}
\label{lin 14 a}\begin{aligned}
\|\check\varphi^\parallel\|_{\mathcal E^{2, \alpha}_\mu(\pD\times \R)}&\leq
C\ve^{{1}-\alpha}\left\{\|(\widetilde \chi
g)^\parallel\|_{\mathcal E^{0, \alpha}_\mu(\pD\times \R)} +\ve^{-1}\delta(\ve)\|(\widetilde\chi\ g)^\perp\|_{\mathcal C^{0, \alpha}_\mu(\pD\times \R)} 
\right.\\
&\qquad\qquad \left.+\ve^{-3-\alpha}\delta(\ve)\|(1-\chi^*)g\|_{\mathcal C^{0, \alpha}(\R^{d-1}\times
S_{2T_\tau})}\right\},\medskip \\
\|\check\varphi^\perp\|_{\mathcal C^{2, \alpha}_\mu(\pD\times \R)}&\leq C\ve^{-1}\left\{\|(\widetilde \chi g)^\perp\|_{\mathcal C^{0, \alpha}_\mu(\pD\times \R)}+\ve^{1-\alpha}|(\widetilde\chi\ g)^\parallel\|_{\mathcal E^{0, \alpha}_\mu(\pD\times \R)}
+\ve^{-2-\alpha}\|(1-\chi^*)g\|_{\mathcal C^{0, \alpha}(\R^{d-1}\times S_{2T_\tau})}\right\},
\end{aligned}
\end{equation}
and
\begin{equation}
\label{lin 14 b}
\|\psi\| _{\mathcal C^{2, \alpha}(\R^{d-1}\times S_{2T_\tau})}\leq C\ve^{-1-\alpha}\|(1-\chi^*)g\|_{\mathcal C^{0, \alpha}(\R^{d-1}\times S_{2T_\tau})}+\mathcal O(e^{\,-c\ve^{-1/8}})\|\widetilde \chi g\|_{\mathcal C^{0, \alpha}_\mu(\pD\times \R)}.
\end{equation}

}}

Using  these {\it a priori} estimates we can solve the system  (\ref{lin 12})--(\ref{lin 13}) by a standard fixed point argument. To do this we replace  the functions $\check\varphi^\parallel, \check\varphi^\perp, \psi$ on the right hand side of the system by known functions $\check\varPhi^\parallel, \check\varPhi^\perp, \varPsi$ which  satisfy estimates of the same type as (\ref{lin 14 a})--(\ref{lin 14 b}) but with constants bigger that those appearing in (\ref{lin 14 a})--(\ref{lin 14 b}). Then we have a map
\[
(\check\varPhi^\parallel, \check\varPhi^\perp, \varPsi)\longmapsto (\check\varphi^\parallel, \check\varphi^\perp, \psi),
\]
from a certain ball in the space $\mathcal E^{2, \alpha}_\mu(\pD\times \R)\times \mathcal C^{2, \alpha}_\mu(\pD\times \R)\times \mathcal C^{2, \alpha}(\R^{d-1}\times S_{2T_\tau})$ into itself. This and the Lipschitz character  of this map being evident from the way we have derived {\it a priori} estimates allows for an  application of the Banach fixed point theorem. We leave the details to the reader and simply state  this as a result for (\ref{lin 11}).
\begin{lemma}\label{exist lin}
For each sufficiently small $\ve$ there exists a solution of (\ref{lin 11}) in the form $\varphi=\chi^*\check\varphi\circ Y_{\ve, h}+\psi$ such that estimates (\ref{lin 14 a})--(\ref{lin 14 b}) hold. 
\end{lemma}


\subsection{Proof of Theorem \ref{teorema 1}}
Now we can finish solving the nonlinear problem 
\begin{equation}
\label{lsr 1}
L_\ve\varphi={\ell}_\ve-N_\ve(w^*)-Q_\ve(\varphi).
\end{equation}
As we saw above we need to modify this equation by introducing Lagrange multipliers. Thus we will consider:
\begin{equation}
L_\ve\varphi=\ell_\ve-N_\ve(w^*)-Q_\ve(\varphi)+\chi^*\sum_{j=1}^d{\tt c}_j Z^{T, {\tt e}_j}_{\tau, \ve}.
\label{lsr 2}
\end{equation}
To solve this problem we use a fixed point argument and the linear theory  in Lemma \ref{exist lin} above. 
The first task is to calculate  the size of the error of the approximation $\ell_\ve-N_\ve(w^*)$. This is straightforward using the definition of $w^*$ and  formula (\ref{fermi 2}). We recall here that $h=\ve^2 h_0$, where $h_0$ is a constant and consequently this last formula simplifies significantly. We can write:
\[
\ell_\ve-N_\ve(w^*)=\chi^*[\ell_\ve-N_\ve(w)]+[N(w^*)-\chi^*N(w)-(1-\chi^*)N({\mathbb H})]
\equiv A_1+A_2,
\]
since $\ell_\ve=N_\ve({\mathbb H})$ in $\mathrm{supp}\,(1-\chi^*)$. Using  exponential  decay of $w-(\pm 1+\sigma_\ve)$ when $\st\to \pm \infty$ we get easily:
\begin{equation}\label{lsr 3}
\begin{aligned}
\|Y_{\ve, h}^*\widetilde \chi A_2\|_{\mathcal C^{0, \alpha}_\mu(\pD\times \R)}&\leq C_0 e^{\,-c_\mu\ve^{-1/3}},\\
\|(1-\chi^*)A_2\|_{\mathcal C^{0, \alpha}(\R^{d-1}\times S_{2T_\tau})}&\leq C_0 e^{\,-\theta\ve^{-1}}.
\end{aligned}
\end{equation}
To estimate $A_1$ some standard calculations which we will omit are needed (c.f Section \ref{sec formal}). As a result we get
{{
\begin{equation}\label{lsr 4}
\begin{aligned}
\|[Y_{\ve, h}^*\widetilde \chi A_1]^\parallel\|_{\mathcal E^{0, \alpha}_\mu(\pD\times \R)}&\leq C_0\ve^3, \\
\|[Y_{\ve, h}^*\widetilde \chi A_1]^\perp\|_{\mathcal C^{0, \alpha}_\mu(\pD\times \R)}&\leq C_0\ve^3,\\
\|(1-\chi^*)A_1\|_{\mathcal C^{0, \alpha}(\R^{d-1}\times S_{2T_\tau})}&\leq C_0 e^{\,-\theta\ve^{-1}},
\end{aligned}
\end{equation}
}}
where $C_0$, $c_\mu$ and $\theta$ are positive constants.
Now we use the linear theory developed in the previous section to solve the nonlinear problem (\ref{lsr 2}). Thus we write $\varphi=\chi^*\check \varphi\circ Y_{\ve, h}+\psi$, and further decompose $\check \varphi=\check \varphi^\parallel+\check \varphi^\perp$ where $\check \varphi^\parallel\in \mathcal X\cap \mathcal C^{2, \alpha}_\mu(\pD\times \R)$, $\varphi^\perp\in \mathcal Y\cap \mathcal C^{2, \alpha}_\mu(\pD\times \R)$ and $\psi\in \mathcal C^{2, \alpha}(\R^{d-1}\times S_{2T_\tau})$. To set up a fixed point scheme we fix  functions $\widetilde\varphi^\parallel$, $\widetilde\varphi^\perp$ and $\widetilde \psi$ in these sets such that 
{{
\begin{equation}\label{lsr 5}
\begin{aligned}
\|\widetilde\varphi^\parallel\|_{\mathcal E^{2, \alpha}_\mu(\pD\times \R)}&\leq K\ve^{4-\alpha}, \\
\|\widetilde\varphi^\perp\|_{\mathcal C^{2, \alpha}_\mu(\pD\times \R)}&\leq K\ve^2,\\
\|\widetilde \psi\|_{\mathcal C^{2, \alpha}(\R^{d-1}\times S_{2T_\tau})}&\leq K e^{\,-\bar \theta\ve^{-1}},
\end{aligned}
\end{equation}
}}
where $K$ is a large constant to be chosen and $\bar\theta\in (\theta/2, \theta)$ is a constant. Let us denote the right hand side of (\ref{lsr 2}) by $\mathfrak g$. It is evident that under the assumptions  (\ref{lsr 5}), and with a suitable choice of the constants $\alpha>0$ and $\mu\in (0, |\eta|)$  we can solve the problem (\ref{lsr 2}) for functions $(\check \varphi^\parallel,\check \varphi^\perp, \psi)$ which again satisfy (\ref{lsr 5}). Thus we have a non-linear map 
\[
(\widetilde\varphi^\parallel, \widetilde\varphi^\perp, \widetilde \psi)\longmapsto (\check \varphi^\parallel,\check \varphi^\perp, \psi),
\]
of this set into itself. To show that this map is a contraction is straightforward, using the quadratic nature of the nonlinear function $Q(\varphi)$. 
At the end we have a solution of the problem:
\begin{equation}
\label{balance 0}
\ve \Delta u+\frac{1}{\ve} f(u)=\ell_\ve + \sum_{j=1}^d \chi^* {\tt c}_j Z^{T, {\tt e}_j}_{\tau,\ve},\quad \mbox{in}\ \R^{d-1}\times S_{2T_\tau}
\end{equation}
where $Z^{T,{\tt e}_j}_{\tau, \ve}$ is the  (approximate) element of the kernel of the linear operator ${\mathbb L}_\ve$  associated with translation in the direction of the $x_j$ axis, see (\ref{ls 3}).  To show that in fact
\[
{\tt c}_j=0, \quad j=1, \dots, d,
\]
we need:
\begin{lemma}[Balancing formula]
Let $X=\sum a_j \partial_{x_j}$ be the infinitesimal generator of translations or rotations in $\R^d$. For any $\mathcal C^2(\R^d)$ function  it holds:
\begin{equation}
\label{balance 1}
\mathrm{div}\, \left[\left(\frac{\ve}{2}|\nabla u|^2-\frac{1}{\ve} F(u)\right) X(u)-\ve X(u)\nabla u\right]=-[\ve\Delta u+\frac{1}{\ve}F'(u)]X(u).
\end{equation}
\end{lemma}
We will take $X_j=\partial_{x_j}$ for some $1\leq j\leq d$ and integrate  the balancing formula  over the cylinder $\mathcal C_{R}=B_R\times S_{2T_\tau}$.  Using (\ref{balance 0}) and Green's theorem  we  get:
\[
\underbrace{\int_{\partial\mathcal C_{R}}\left(\frac{\ve}2|\nabla u|^2-\frac{1}{\ve}F(u)+\ell_\ve u\right)n_j\, dS}_{I_R}-\underbrace{\int_{\partial\mathcal C_{R}}\partial_{x_j} u \partial_n u\,dS}_{II_R}= - \underbrace{\int_{\mathcal C_{R}}\left(\sum_{j'=1}^d\chi^* {\tt c}_{j'} Z^{T, {\tt e}_{j'}}_{\tau,\ve}\right)\partial_{x_j} u}_{{III}_R}.
\]
The first integral $I_R$ is $0$ on the top and the bottom of $\mathcal C_R$ and on the other hand, using the asymptotic behavior of the solution we get finally $\lim_{R\to \infty}I_R=0$. In the second integral the integrals over the top and the bottom of $\mathcal C_R$ cancel because $u$ is periodic. Then, from exponential decay of the derivatives of $u$ we get $\lim_{R\to \infty}II_R=0$. Finally, we note that 
\[
\partial_{x_j} u\approx Z^{T, {\tt e}_j}_{\tau, \ve},
\]
hence 
\[
III_R={\tt c}_j \int_{\pD\times \R}|Z^{T, {\tt e}_j}_{\tau, \ve}|^2+o(1)\sum_{j'=1}^d{\tt c}_{j'},
\]
from which we get immediately ${\tt c}_j=0$, $j=1, \dots, d$. 



\end{document}